\documentclass{article}
\usepackage{amsfonts}
\usepackage{amsmath}
\usepackage{graphicx}
\usepackage{bibspacing}

\usepackage{amsfonts}
\usepackage{amssymb,amsmath,amsthm, amsfonts}
\usepackage{color}
\usepackage[dvipsnames]{xcolor}
\def\sqr#1#2{{\vcenter{\vbox{\hrule height.#2pt
              \hbox{\vrule width.#2pt height#1pt \kern#1pt \vrule width.#2pt}
          \hrule height.#2pt}}}}

\def\esssup{\mathop{\rm esssup}}
\def\essinf{\mathop{\rm essinf}}

\def\exp{\mathop{\rm exp}}
\def\inf{\mathop{\rm inf}}

\font\tenbb=msbm10 \font\sevenbb=msbm7 \font\fivebb=msbm5
\newfam\bbfam
\scriptscriptfont\bbfam=\fivebb \textfont\bbfam=\tenbb
\scriptfont\bbfam=\sevenbb

\newtheorem{lemma}{Lemma}[section]
\newtheorem{remark}{Remark}[section]
\newtheorem{example}{Example}[section]
\newtheorem{theorem}{Theorem}[section]
\newtheorem{corollary}{Corollary}[section]

\newtheorem{definition}{Definition}[section]
\newtheorem{proposition}{Proposition}[section]

\makeatletter
   
   \@addtoreset{equation}{section}
\makeatother

\begin{document}

\title{Backward stochastic differential equations with conditional reflection and related recursive optimal control problems}

\author{
Ying Hu \thanks{Institut de Recherche Math\'ematique de Rennes, Universit\'e Rennes 1,
 35042 Rennes Cedex, France.  Email:{ying.hu@univ-rennes1.fr}. This author's research is partially supported by Lebesgue Center of Mathematics ``Investissements d'avenir" Program (No. ANR-11-LABX-0020-01), by ANR CAESARS
(No. ANR-15-CE05-0024) and by ANR MFG (No. ANR-16-CE40-0015-01)}
 \and
Jianhui Huang \thanks{ Department of Applied Mathematics, The Hong Kong Polytechnic University,
 Hong Kong, China.  Email: {james.huang@polyu.edu.hk}. This author's research is partially supported by RGC Grant PolyU 15301119, 15307621, N PolyU504/19, NSFC 12171407 and KKZT.}
 \and
Wenqiang Li\thanks{Corresponding author. School of Mathematics and Information Sciences, Yantai University,
  264005 Yantai, P. R. China.  Email: {wqliytu@163.com}. This author's research is partially supported by the NSF of P.R. China (No. 12101537) and Doctoral Scientific Research Fund of Yantai University (No. SX17B09).}
}

\date{\today}
\maketitle

\begin{abstract}We introduce a new type of reflected backward stochastic differential equations (BSDEs) for which the reflection constraint is imposed on its main solution component, denoted as $Y$ by convention, but in terms of its conditional expectation $\mathbb{E}[Y_t|\mathcal{G}_{t}]$ on a general sub-filtration $\{\mathcal{G}_{t}\}.$  We thus term such equation as \emph{conditionally reflected BSDE} (for short, conditional RBSDE). Conditional RBSDE subsumes classical RBSDE with a pointwise reflection barrier, and the recent developed BSDE with a mean reflection constraint, as its two special and extreme cases: they exactly correspond to $\{\mathcal{G}_{t}\}$ being the full filtration to represent complete information, and the degenerated filtration to deterministic scenario, respectively. For conditional RBSDE, we obtain its existence and uniqueness under mild conditions by combining the Snell envelope method with  Skorokhod lemma. We also discuss its connection, in the case of linear driver, to a class of optimal stopping problems in presence of partial information. As a by-product, a new version of comparison theorem is obtained. With the help of this connection, we study \emph{weak formulations} of a class of optimal control problems with reflected recursive functionals by characterizing the related optimal solution and value. Moreover, in the special case of recursive functionals being RBSDE with pointwise reflections, we study the \emph{strong formulations} of related stochastic backward recursive control and zero-sum games, both in non-Markovian framework, that are of their own interests and have not been fully explored by existing literatures yet.
\end{abstract}




\section{Introduction}

Reflected backward stochastic differential equations (RBSDEs) were firstly introduced by El Karoui, Kapoudjian, Pardoux, Peng and Quenez in \cite{KKPPQ} for which the solution is an adapted triple processes $(Y,Z,K)$ satisfying the following backward stochastic system in an integral form
\begin{equation}\label{BSDE}
 Y_t=\xi+\int_t^Tf(s,Y_s,Z_s)ds+K_T-K_t-\int_t^TZ_sdW_s,\ 0\leq t\leq T,\\
\end{equation}
subject to a \emph{pointwise} constraint\begin{equation}\label{constraint-peng}
 Y_t\geq S_t,\ 0\leq t\leq T,
\end{equation}for a given barrier process $S$. The term $K$ of the solution is used to push the \emph{main} solution component $Y$ to satisfy the constraint \eqref{constraint-peng} in a minimum energy way, i.e., $\int_0^T(Y_t-S_t)dK_t=0$. In \cite{KKPPQ}, the authors show the well-posedness of the solution $(Y,Z,K)$ of the above equation for a given terminal condition $\xi$, a Lipschitz driver $f,$ and a continuous barrier process $S$. They also establish its connection to both optimal stopping problems and related obstacle problems of parabolic partial differential equations. Due to its interesting structure, RBSDEs have been extensively applied, among others, into the problems such as pricing of the American option \cite{EPQ-1}, and dynamic recursive portfolio problems \cite{EPQ-2,WY2008}.

Recently, Briand, Elie and Hu \cite{BEH2018} introduced BSDEs with mean reflection, which is a  type of RBSDEs satisfying \eqref{BSDE} but subject to a constraint condition in terms of the expectation as
\begin{equation}\label{constraint-Hu}
 \mathbb{E}[\ell(t,Y_t)]\geq 0,\ 0\leq t\leq T,
\end{equation}for some given loss function $\ell$.
In contrast to the pointwise reflection constraint \eqref{constraint-peng}, condition \eqref{constraint-Hu} is described in sense of the distribution of the term $Y_t$ at each instant time $t$. In \cite{BEH2018}, the authors construct a unique solution $(Y,Z,K)$ with requiring $K$ to be deterministic under appropriate conditions on the data $(\xi,f,\ell)$. Using BSDEs with the above mean reflection, they studied the related super-hedging of a claim under a given running static risk management constraint. Since then, many extension works on BSDEs with mean reflection have been studied, among others, \cite{HHLLW,HMW2022} with quadratic growth in $z$ of the driver $f$, \cite{LW2019} with BSDEs driven by $G$-Brownian motion, and \cite{BH2021} with  related particles systems.

In this paper, we study a new type of RBSDEs, which is called \emph{conditional RBSDEs}, for which the reflection barrier is defined via a general conditional expectation operator, that is, system \eqref{BSDE} is subject to a constraint of the form:
\begin{equation}\label{constraint}
\mathbb{E}[Y_t-S_t|\mathcal{G}_t]\geq 0,\ 0\leq t\leq T,
\end{equation}
for some generic $\mathcal{G}_t\subseteq\mathcal{F}_t, 0\leq t\leq T$, so the sub-filtration $\mathbb{G}=\{\mathcal{G}_t\}_{0\leq t\leq T}$ (see Section 2 for more details) stands for \emph{partial information} that is common for various real decision making applications. It is worth to note that both the constraints \eqref{constraint-peng} and \eqref{constraint-Hu} (with linear loss function $\ell$) can be seen as the special case of condition \eqref{constraint}. In fact, reflection condition \eqref{constraint-peng} (resp., \eqref{constraint-Hu}) corresponds to the \emph{full information} (resp., \emph{(degenerated) deterministic  scenario}) situation in \eqref{constraint} when $\mathcal{G}_t=\mathcal{F}_t$ (resp., $\mathcal{G}_t=\mathcal{F}_0$), $0\leq t\leq T$.

Condition \eqref{constraint} is strongly suggested by portfolio selection problems subject to some state constraints but in the context of partial information. On one hand, notice that the partial-information feature in portfolio investments have been extensively studied by many mathematical finance works such as  \cite{BDL2010,HY2006,L1995,XZ2007,ZXZ2020}. Roughly speaking, in these studies, investors are often posed in a situation where only part of the overall information of the market can be accessed;
for instance, the driving noise information available to the investor is often incomplete due to some observation or measurement limitations, so the investors can only observe a subset of underlying noise components but not all. This is typical for various commonly-seen real situations, especially when some latent factors are indispensable to drive the dynamic evolution but cannot be accurately and instantaneously calibrated. On the other hand, due to some regulation criteria or behavior pattern, investors (e.g., fund managers) should make sure the state (fund account) to be controlled or steered to meet some constraint (e.g., above some market benchmark or average level) for the purposes such as principal evaluation or stimulus mechanism. This brings some obstacle constraint on the underlying state evolution. Together, some constraint portfolio selections with partial information are hence suggested. For illustration, we present two examples below to formulate our RBSDE \eqref{BSDE} with constraint \eqref{constraint}, and motivate related optimal control problems.

\emph{Example 1.1}. (Pricing American options with partial information) Denote ${\mathbb{G}}$ as a generic \emph{partial information} a representative agent can access from the market information sources. The pricing to an American contingent claim at each time $t$, consists of the selection of a stopping time $\tau$ and a payoff $\widetilde{S}_\tau$ on exercise. It is natural to restrict that both the stopping time $\tau$ and stopped (truncated) process $\widetilde{S}=\{\widetilde{S}_\tau\}$ to be ${\mathbb{G}}$-adapted because they are both chosen upon the information available to a specific decision maker (i.e., agent) from the market. We denote by $S$ the payoff of the American claim for an \emph{idealized} agent who may access the full information of the market. Then, for a \emph{realistic} agent with only partial information, it is reasonable and more practical to only anticipate the conditional, unbiased payoff: $\widetilde{S}_t=\mathbb{E}[S_t|\mathcal{G}_t]$, $0\leq t\leq T,$ in current information capacity.

It is well known that for each selection on $\tau$, there exists a unique strategy $(y^\tau,z^\tau)$ to replicate the payoff ${S}$, where $(y^\tau,z^\tau)$ is the unique solution of the BSDE
$$-dy^\tau_t=f(t,y^\tau_t,z^\tau_t)dt-z^\tau_t dW_t,\ y^\tau_\tau={S}_\tau,$$
for some convex and Lipschitz driver $f$. Then the optimal pricing of the American contingent claim $\widetilde{S}$ with partial information at each time $t$, is given by
$$y_t:=\esssup_{\tau\in[t,T]:{\mathbb{G}}\text{-adapted}}
\mathbb{E}[y_t^\tau|\mathcal{G}_t].$$
Similar to the full information context studied in \cite{KKPPQ}, we may expect (in fact, we verify the case when the driver $f$ is linear in $(y,z)$ in Section 3) that the value process $y$ can be characterized as follows
$$y_t=\mathbb{E}[Y_t|\mathcal{G}_t],\ 0\leq t\leq T,$$
where $(Y,Z,K)$ is the solution triple of conditional RBSDE \eqref{BSDE} along with \eqref{constraint}.
Condition \eqref{constraint} means that the term $Y$ is always required to be larger than the payoff $\widetilde{S}$ at each time $t$, in terms of conditional expectation on the available information $\{\mathcal{G}_t\}$. Then we can use conditional RBSDE \eqref{BSDE} and \eqref{constraint} to price American options with partial information. Indeed, when all agents are symmetric in their information  (i.e., all access the same sub-filtration $\{\mathcal{G}_t\}$), some equilibrium on supple-demand condition will be achieved by game-theoretic analysis, and the market price will be thus formalized on such partial information basis.
Specifically, in one extreme and idealized case when agents can all access full information, the above result will reduce to the pricing results studied by El Karoui, Pardoux and Quenez \cite{EPQ-1}.
Another extreme case is when all agents only access trivial filtration (i.e., cannot observe any realization of stochastic scenarios), then an optimal (deterministic) stopping time arises to get an \emph{expected} payoff evaluated by the agent, see Example \ref{eg4.1}.



\emph{Example 1.2}. (Recursive reflected utility maximization with partial information)
We consider an optimal portfolio selection problem in a market consisting of a risk-free bond and $d$ risky assets. The bond is assumed to be zero interest rate and the discounted (by the bond) individual risky asset price $V_t^i$ at time $t$ has the following form
\begin{equation}\label{101601}
\frac{dV_t^i}{V_t^i}=b_t^idt+\sigma_t^idW_t^i,\ i=1,2,\cdots, d.
\end{equation}
Here, $b^i,\sigma^i,i=1,2,\cdots,d,$ are given return and volatility rates, respectively, and $W=(W^1,\cdots,W^d)$ is a standard $d$-dimensional Brownian motion. An investor only observes the (public market) price of listed risky assets (e.g., stocks) $1,2,\cdots,m$ (with $1\leq m\leq d$) since, for instance, the prices of some unlisted risky assets are latent as described in principle-agent situation with partial information or hidden actions by Williams \cite{W2009}.  Then the partial information available to the investor in public market is$$\mathcal{G}_t=\sigma\{V_s^i,0\leq s\leq t,\ i=1,2,\cdots,m\},\ 0\leq t\leq T,$$
which is same to $\sigma\{W_s^i,0\leq s\leq t,\ i=1,2,\cdots,m\}$ when both $b^i$ and $\sigma^i$, $i=1,2,\cdots,m,$ are deterministic.
Let $\pi_t=(\pi_t^1,\cdots,\pi_t^d)$ be the proportion of the amount invested in risky assets at time $t$, which should be ${\mathbb{G}}$-adapted. Then the wealth process $X_t^\pi$ at time $t$ with the proportion $\pi$ should satisfy$$dX_t^\pi=\sum_{i=1}^d\pi_t^iX_t^\pi[b_t^idt+\sigma_t^idW_t^i].$$
Such portfolio model with partial information has been studied by Nagai and Peng \cite{NP2002} when addressing a type of risk-sensitive optimization problems on an infinite time horizon.

The expected utility of the investor is of recursive utility, denoted by $Y^\pi_0$, which can be described by the solution of classical BSDE according to \cite{KPQ1997}. Moreover, let the constant $a$ be the minimum utility threshold acceptable to the investor, i.e., $Y^\pi_0\geq a,$ for all admissible investment proportion $\pi$. In fact, we can consider a stronger dynamic constraint as follows, which depends on the evaluation of the utility at each time $t$ based on the partial information $\mathbb{G}$, i.e.,
$$\mathbb{E}[Y^\pi_t|\mathcal{G}_t]\geq \widetilde{S}_t,\ 0\leq t\leq T,$$
where $\widetilde{S}$ is $\mathbb{G}$-adapted with $\widetilde{S}_0=a$. Combining the above two factors, we obtain a recursive reflected utility $Y^\pi_0$, which is defined by the solution of the following controlled conditional RBSDE
\begin{equation}\nonumber
\left\{
\begin{array}{ll}
\displaystyle Y_t^\pi=\Phi(X_T^\pi)+\int_t^Tf(s,X_s^\pi,Y_s^\pi,Z_s^\pi)ds+K_T^\pi-K_t^\pi-\int_t^TZ_s^\pi dW_s,\ 0\leq t\leq T,\\
\displaystyle \mathbb{E}[Y^\pi_t|\mathcal{G}_t]\geq \widetilde{S}_t,\ 0\leq t\leq T,\\
 \end{array}
 \right.
\end{equation}
where $\Phi$ and $f$ represents the terminal and running utility, respectively.
The aim of the investor is to choose an admissible proportion $\pi$ to maximize the recursive $Y_0^\pi$ using the available information $\mathbb{G}$ only.

Inspired by the above examples, we aim to study the well-posedness of conditional RBSDE consisting of \eqref{BSDE}, \eqref{constraint}, and associated backward recursive reflected control problems with partial information. We also address the counterpart of control problems in full information case which may admit more explicit results under more relaxed assumptions. In order to guarantee the uniqueness of the solution $(Y,Z,K)$, we consider the case that $K$ is required to be $\mathbb{G}$-adapted as explained in Remark \ref{re2.1}. Similar to classical RBSDE studied in \cite{KKPPQ} and BSDEs with mean reflection in \cite{BEH2018}, the form of conditional RBSDEs involves in a flat condition: $\int_0^T\mathbb{E}[Y_t-S_t|\mathcal{G}_t]dK_t=0$ besides \eqref{BSDE} and \eqref{constraint}. It is worth to point out that, when studying the well-posedness of the solution, the partial information is only required to satisfy the usual filtration condition and the additional left-quasi-continuous condition (to ensure the conditional expectation $\mathbb{E}[\cdot|\mathcal{G}_t]$ is continuous in $t$).
For related backward recursive reflected control problems with partial information, we consider the partial information constructed by a subset of components of driving Brownian motion noises, which is motivated by \emph{Example 1.2.} We emphasize that this construction approach includes a large class of partial information models, as explained in Remark \ref{re4.1}.

The rest of this paper is organized as follows. Section 2 formulates the conditional RBSDE on a general sub-filtration along with necessary assumptions. We then study the well-posedness (including the existence and the uniqueness as well as a prior estimate) of the solution of conditional RBSDE. Section 3 is devoted to the connection between conditional RBSDE and a new class of optimal stopping problems in presence of partial information. As a byproduct, a related comparison theorem is also derived that has its own interests in theoretical analysis. Section 4 considers two types of backward recursive reflected control problems with partial information in case of the driver for recursive functional being linear and convex, respectively. Section 5 continues to study the strong formulations of nonlinear backward recursive reflected functionals for which both stochastic control and zero-sum game problems in non-Markovian framework are examined. Some equivalence between strong and weak formulations is also established.

\section{Conditionally reflected BSDEs}
\subsection{Preliminary}

Let $T>0$ be a finite time horizon.  Suppose that $\{W_s,\ s\in[0,T]\}$ is a $d$-dimensional standard Brownian motion defined on a  probability space
$(\Omega,{\mathcal{F}},\mathbb{P})$. We denote by $\mathbb{E}$ the (conditional) expectation under the probability measure $\mathbb{P}$ and by $\mathbb{F}=(\mathcal{F}_s)_{s\geq 0}$ the complete filtration generated by Brownian motion $W$. Let
$$\mathcal{G}_t\subseteq\mathcal{F}_t,\ t\in[0,T],$$
be a given sub-filtration of $\mathcal{F}_t$ satisfying the following \emph{basic assumption}\footnote{
The condition (i) is essential and classical for the information filtration and the condition (ii) is introduced to ensure the continuity property of conditional expectation $\mathbb{E}[\cdot|\mathcal{G}_t]$  in $t$.}:\\
(i) the usual condition (i.e., non-decreasing and right-continuous);\\
(ii) left-quasi-continuous (i.e., left-continuous with respect to stopping times).\\
We denote this sub-filtration by $\mathbb{G}=(\mathcal{G}_t)_{0\leq t\leq T}$ and refer it as \emph{partial information} (resp., $\mathbb{F}$ as \emph{full information}) inspired by  examples in Introduction.  The trivial $\sigma$-field is denoted by $\mathcal{H}$, i.e.,
$\mathcal{H}=\mathcal{F}_0=\mathcal{G}_0$, which is referred as \emph{(degenerated) deterministic scenario}.
 We  introduce   the following spaces of processes:
\begin{itemize}

\item ${\mathcal{S}}^{2}(0,T;\mathbb{R})=\Big\{\varphi \Big| \varphi :\Omega\times [0,T]
\rightarrow \mathbb{R}$ is  $\mathbb{F}$-adapted and continuous: $\|\varphi\|_{\mathcal{S}^2}^2=\mathbb{E}[\sup\limits_{t\in [0,T]}|\varphi_{t}|^{2}]<+\infty
\Big\}. $
 \item $\mathcal{H}^{2}(0,T;\mathbb{R}^d)=\Big\{\varphi\Big|\varphi :\Omega\times[0,T]  \rightarrow \mathbb{R}^d$ is $\mathbb{F}$-predictable:
$\|\varphi\|_{\mathcal{H}^2}^2=\mathbb{E}\Big[\displaystyle\int_{0}^{T}|\varphi _{t}|^{2}dt\Big]<+\infty \Big\}.$

\item ${\mathcal{A}}^{2}(0,T;\mathbb{R})=\Big\{\varphi \Big| \varphi\in\mathcal{S}^2(0,T;\mathbb{R})
$, $\varphi$ is nondecreasing, $\varphi_0=0\Big\}. $

\item ${\mathcal{A}}^{2}_{\mathbb{G}}(0,T;\mathbb{R})=\Big\{\varphi \Big| \varphi\in\mathcal{A}^2(0,T;\mathbb{R}),\varphi$ is  $\mathbb{G}$-adapted$\Big\}$.
\end{itemize}
For all these spaces, we write $\mathcal{S}^2,\mathcal{H}^2,\mathcal{A}^2,\mathcal{A}_\mathbb{G}^2$ when there is no confusion hereafter.
We are interested in the following \emph{conditional RBSDE} associated with parameters: the terminal condition $\xi$, the driver $f$ and a barrier process $S$:
\begin{equation}\label{constrained-bsde}
\left\{
\begin{aligned}
&Y_t=\xi+\int_t^T f(s,Y_s,Z_s)ds +K_T-K_t-\int_t^TZ_sdW_s,\ t\in[0,T],\\
&\mathbb{E}[Y_t-S_t|\mathcal{G}_t]\geq 0,\ {\mathbb{P}\text{-a.s.}},\ 0\leq t\leq T;\
 \int_0^T\mathbb{E}[Y_t-S_t|\mathcal{G}_t]dK_t=0,\ { \mathbb{P}\text{-a.s.}}
\end{aligned}
\right.
\end{equation}
It is clear that conditional RBSDE   \eqref{constrained-bsde} will reduce to the classical RBSDE introduced by
El Karoui, Kapoudjian, Pardoux, Peng and Quenez \cite{KKPPQ} when $\mathcal{G}_t=\mathcal{F}_t$, and to BSDE with linear mean reflection studied recently by Briand, Elie and Hu \cite{BEH2018} when $\mathcal{G}_t\equiv\mathcal{H}$, $0\leq t\leq T$.

\begin{definition}
A solution of conditional RBSDE  \eqref{constrained-bsde}
is a triple of processes $(Y,Z,K)\in\mathcal{S}^2\times
\mathcal{H}^2\times\mathcal{A}_\mathbb{G}^2$ satisfying \eqref{constrained-bsde}.
\end{definition}
\begin{remark}\label{re2.1}
As  the example displayed in the Introduction of \cite{BEH2018} shows that, we can not expect to obtain the uniqueness of the solution if we allow $K\in\mathcal{A}^2$. As a result, we restrict ourself to find the term $K$ of the solution  in the space $\mathcal{A}_\mathbb{G}^2$ instead of $\mathcal{A}^2$.
\end{remark}

We introduce the basic assumptions of parameters $(\xi,f,S)$ of conditional RBSDE. Let the mapping
$$
f:\Omega\times[0,T]\times\mathbb{R}\times\mathbb{R}^d\rightarrow\mathbb{R},\
 $$
be $\mathcal{P}\otimes\mathcal{B}(\mathbb{R})\otimes\mathcal{B}(\mathbb{R}^d)$
-measurable, where $\mathcal{P}$ stands for the $\sigma$-algebra of $\mathbb{F}$-progressive subsets of $\Omega\times [0,T]$. Suppose that these parameters satisfy
\begin{equation}\nonumber
{\bf (H1)}
\left\{
\begin{array}{l}
\text{(i)\  There\ exists\ a\ constant}\
\mu>0\ \text{such\ that,\ for\ all}\ (t,\omega)\in[0,T]\times\Omega,\ (y,z), (y',z')\in\mathbb{R}^{1+d},\\
\qquad  |f(t,y,z)-f(t,y',z')|\leq \mu(|y-y'|+|z-z'|),\ \mathbb{P}\text{-a.s.,}\
\text{and}\ \mathbb{E}[\int_0^T|f(t,0,0)|^2dt]<\infty.\\
\text{(ii)\  The\ barrier\ process}\
S\ \text{is\ in}\ \mathcal{S}^2. \\
\text{(iii)}\ \text{The\ terminal\ condition}\ \xi\in L^2(\Omega,\mathcal{F}_T,\mathbb{P})\ \text{ such that}\
 \mathbb{E}[\xi-S_T|\mathcal{G}_T]\geq 0,\ \mathbb{P}\text{-a.s.}
 \end{array}
 \right.
 \end{equation}
The remainder of this section is devoted to the  study of the well-posedness of the solution of conditional RBSDE \eqref{constrained-bsde} under the Assumption (H1).

\subsection{Uniqueness of the solution}

Since the term $K$ of the solution is required to be $\mathbb{G}$-adapted (see Remark \ref{re2.1}), we first derive its expression in terms of the conditional expectation with respect to partial information $\mathbb{G}$ by using Skorohod lemma.
\begin{proposition}\label{prop3.2}
Suppose that $(Y,Z,K)\in\mathcal{S}^2\times\mathcal{H}^2\times\mathcal{A}^2_\mathbb{G}$ is a solution of conditional RBSDE \eqref{constrained-bsde}. Then the term $K$ has the following representation: for $t\in[0,T]$ and each $\omega\in\Omega$,
\begin{equation}\label{091501}
\begin{aligned}
(K_T-K_t)(\omega)
=&\sup_{t\leq s\leq T}\Big(\mathbb{E}[\xi|\mathcal{G}_T]+
\mathbb{E}[\int_0^Tf(r,Y_r,Z_r)dr|\mathcal{G}_T]-\mathbb{E}[\int_0^{s}
f(r,Y_r,Z_r)dr|\mathcal{G}_{s}]\\
&\qquad\quad+\mathbb{E}[\int_0^{s}Z_rdW_r|\mathcal{G}_{s}]
-\mathbb{E}[\int_0^TZ_rdW_r|\mathcal{G}_T]
-\mathbb{E}[S_{s}|\mathcal{G}_{s}]
\Big)^-(\omega).
\end{aligned}
\end{equation}
\end{proposition}
\begin{proof}
From \eqref{constrained-bsde}, we get
\begin{equation}\nonumber
\begin{aligned}
\mathbb{E}[Y_t|\mathcal{G}_t]
=&\Big(Y_0-\mathbb{E}[\int_0^tf(s,Y_s,Z_s)ds|\mathcal{G}_t]
+\mathbb{E}[\int_0^tZ_sdW_s|\mathcal{G}_t]\Big)-K_t,\ 0\leq t\leq T,\\
\end{aligned}
\end{equation}
which implies that
\begin{equation}\label{091502}
\begin{aligned}
\mathbb{E}[Y_t|\mathcal{G}_t]
-\mathbb{E}[\xi|\mathcal{G}_T]
=&\Big(\mathbb{E}[\int_0^Tf(s,Y_s,Z_s)ds|\mathcal{G}_T]
-\mathbb{E}[\int_0^tf(s,Y_s,Z_s)ds|\mathcal{G}_t]\\
&+\mathbb{E}[\int_0^tZ_sdW_s|\mathcal{G}_t]-
\mathbb{E}[\int_0^TZ_sdW_s|\mathcal{G}_T]\Big)+K_T-K_t.
\end{aligned}
\end{equation}
By putting
\begin{equation}\nonumber
\begin{aligned}
x_t=&\Big(\mathbb{E}[\xi|\mathcal{G}_T]+
\mathbb{E}[\int_0^Tf(s,Y_s,Z_s)ds|\mathcal{G}_T]-\mathbb{E}[\int_0^{T-t}
f(s,Y_s,Z_s)ds|\mathcal{G}_{T-t}]+\mathbb{E}[\int_0^{T-t}Z_sdW_s|\mathcal{G}_{T-t}]\\
&
-\mathbb{E}[\int_0^TZ_sdW_s|\mathcal{G}_T]
-\mathbb{E}[S_{T-t}|\mathcal{G}_{T-t}]
\Big)(\omega),\
y_t=\mathbb{E}[Y_{T-t}-S_{T-t}|\mathcal{G}_{T-t}](\omega),\ k_t=(K_T-K_{T-t})(\omega),
\end{aligned}
\end{equation}
we have from \eqref{091502} that
$
y_t
=x_t
+k_t,\ t\in[0,T].
$
Moreover, the reflection and flat conditions in \eqref{constrained-bsde} mean that $y_t\geq 0$, $\int_0^Ty_tdk_t=0$.
Note that $x_t$ is continuous with respect to $t\in[0,T]$ and $x_0\geq 0$, from Skorohod Lemma, we get  $k_t=\sup_{0\leq s\leq t}x_s^-$, i.e., \eqref{091501}.
\end{proof}
\begin{remark}
When the available information $\mathbb{G}$ is chosen to be full, i.e., $\mathcal{G}_t=\mathcal{F}_t$, $0\leq t\leq T$, then the representation \eqref{091501} will reduce to  Proposition 2.2 in \cite{KKPPQ}. When the available information $\mathbb{G}$ is chosen to be deterministic scenario, i.e.,  $\mathcal{G}_t=\mathcal{H}$, $0\leq t\leq T$, the expression \eqref{091501} has been used to construct the solution of BSDEs with linear mean reflection  in \cite{BEH2018} (see Subsection 3.3 therein).
\end{remark}
With the help of Proposition \ref{prop3.2}, we  get the following a priori estimate of the solution.
\begin{theorem}\label{le4.1}
For $i=1,2,$
let $(Y^i,Z^i,K^i)\in\mathcal{S}^2\times\mathcal{H}^2\times\mathcal{A}_{\mathbb{G}}^2$ be a solution
of conditional RBSDE \eqref{constrained-bsde} associated with parameters $(\xi^i,f^i,S^i)$ satisfying the Assumption (H1). Then there exists a constant $C$ only depending on $T$ and $\mu$ such that, for any $t\in[0,T]$,
\begin{equation}\nonumber
\begin{aligned}
&\mathbb{E}\Big[\sup_{0\leq s\leq T}|Y_s^1-Y_s^2|^2+\int_0^T
|Z_s^1-Z_s^2|^2ds+\sup_{0\leq s\leq T}|(K_T^1-K_s^1)-(K_T^2-K_s^2)|^2\Big]\\
\leq&
C\mathbb{E}\Big[|\xi^1-\xi^2|^2+\int_0^T|f^1(s,Y_s^2,Z_s^2)-f^2(s,Y_s^2,Z_s^2)|^2ds
+\sup_{0\leq s\leq T}|S_s^1-S_s^2|^2\Big].
\end{aligned}
\end{equation}
\end{theorem}
\begin{proof}
For simplicity of the notations, we denote
\begin{equation}\nonumber
\begin{aligned}
 \Delta f(s)&=f^1(s,Y_s^2,Z_s^2)-f^2(s,Y_s^2,Z_s^2),\ \Delta L=L^1-L^2,\ L=Y,Z,K,\xi,S.
\end{aligned}
\end{equation}
Step 1. We show that
\begin{equation}\label{091603}
\begin{aligned}
&\mathbb{E}\Big[|\Delta Y_t|^2+\int_t^T|\Delta Y_s|^2+
|\Delta Z_s|^2ds+ |\Delta K_T-\Delta K_t|^2\Big|\mathcal{G}_t\Big]\\
\leq&
C\mathbb{E}\Big[|\Delta\xi|^2+\int_t^T|\Delta f(s)|^2ds
+\sup_{t\leq s\leq T}|\Delta S_s|^2\Big|\mathcal{G}_t\Big],\ t\in[0,T],\ \mathbb{P}\text{-a.s.}
\end{aligned}
\end{equation}
For any $\beta>0$, applying It\^o's formula to $e^{\beta t}|\Delta Y_t|^2$ we get
\begin{equation}\label{050805}
\begin{aligned}
&e^{\beta t} |\Delta Y_t|^2+\int_t^Te^{\beta s}(\beta  |\Delta Y_s|^2+|\Delta Z_s|^2)ds\\
\leq &e^{\beta T}|\Delta\xi|^2
+\int_t^Te^{\beta s} |\Delta f(s)|^2ds+\frac18\int_t^Te^{\beta s} |\Delta Z_s|^2ds
+(1+2\mu+8\mu^2)\int_t^Te^{\beta s} |\Delta Y_s|^2ds\\
&+2\int_t^Te^{\beta s} \Delta Y_sd\Delta K_s-2\int_t^Te^{\beta s} \Delta Y_s\Delta Z_sdW_s.
\end{aligned}
\end{equation}
Choosing $\beta=2+2\mu+8\mu^2$, it follows from \eqref{050805} that
\begin{equation}\label{050806}
\begin{aligned}
&e^{\beta t} \mathbb{E}[|\Delta Y_t|^2|\mathcal{G}_t]+
\mathbb{E}[\int_t^Te^{\beta s}|\Delta Y_s|^2ds|\mathcal{G}_t]+\frac78\mathbb{E}[\int_t^Te^{\beta s}|\Delta Z_s|^2ds|\mathcal{G}_t]\\
\leq &\mathbb{E}[e^{\beta T}|\Delta \xi|^2|\mathcal{G}_t]
+\mathbb{E}[\int_t^Te^{\beta s} |\Delta f(s)|^2ds|\mathcal{G}_t]
+2\mathbb{E}[\int_t^Te^{\beta s}\Delta Y_sd\Delta K_s|\mathcal{G}_t]\\
\leq &\mathbb{E}[e^{\beta T}|\Delta\xi|^2|\mathcal{G}_t]
+\mathbb{E}[\int_t^Te^{\beta s} |\Delta f(s)|^2ds|\mathcal{G}_t]
+ e^{\beta T} \Big(\frac{1}{\varepsilon}
\mathbb{E}[\sup_{s\in[t,T]}|\Delta S_s|^2|\mathcal{G}_t]
+\varepsilon\mathbb{E}[|\Delta K_T-\Delta K_t|^2|\mathcal{G}_t]\Big),
\end{aligned}
\end{equation}
where the last inequality follows from
\begin{equation}\label{110401}
\begin{aligned}
&2\mathbb{E}[\int_t^Te^{\beta s}\Delta Y_sd\Delta K_s|\mathcal{G}_t]=2\mathbb{E}[\int_t^Te^{\beta s} (\Delta Y_s-\Delta S_s)d\Delta K_s|\mathcal{G}_t]+
2\mathbb{E}[\int_t^Te^{\beta s} \Delta S_sd\Delta K_s|\mathcal{G}_t]\\
=&2\mathbb{E}[\int_t^Te^{\beta s} {\mathbb{E}[\Delta Y_s-\Delta S_s|\mathcal{G}_s]}d\Delta K_s|\mathcal{G}_t]+
2\mathbb{E}[\int_t^Te^{\beta s}\Delta  S_sd\Delta K_s|\mathcal{G}_t]\\ 
\leq&
2\mathbb{E}[\int_t^Te^{\beta s}\Delta  S_sd\Delta K_s|\mathcal{G}_t]
\leq e^{\beta T} \Big(\frac{1}{\varepsilon}
\mathbb{E}[\sup_{s\in[t,T]}|\Delta S_s|^2|\mathcal{G}_t]
+\varepsilon\mathbb{E}[|\Delta K_T-\Delta K_t|^2|\mathcal{G}_t]\Big).
\end{aligned}
\end{equation}
Since $\Delta K_T-\Delta K_t=\Delta Y_t-\Delta \xi-\int_t^Tf^1(s,Y_s^1,Z_s^1)-f^1(s,Y_s^2,Z_s^2)+\Delta f(s)ds+\int_t^T\Delta Z_sdW_s$, we get
\begin{equation}\label{050901}
\begin{aligned}
&\mathbb{E}[|\Delta K_T-\Delta K_t|^2|\mathcal{G}_t]\\
\leq& C(T,\mu)\Big(\mathbb{E}[|\Delta Y_t|^2|\mathcal{G}_t]
+\mathbb{E}[|\Delta \xi|^2|\mathcal{G}_t]
+\mathbb{E}[\int_t^T |\Delta Y_s|^2+|\Delta Z_s|^2ds|\mathcal{G}_t]
+\mathbb{E}[\int_t^T |\Delta f(s)|^2ds|\mathcal{G}_t]
\Big).
\end{aligned}
\end{equation}
Substituting \eqref{050901} into \eqref{050806}, choosing $\varepsilon$ small enough, we get
\begin{equation}\nonumber
\begin{aligned}
&\mathbb{E}[ |\Delta Y_t|^2|\mathcal{G}_t]+\mathbb{E}[\int_t^T|\Delta Y_s|^2+|\Delta Z_s|^2ds|\mathcal{G}_t]
+\mathbb{E}[|\Delta K_T-\Delta K_t|^2|\mathcal{G}_t]\\
\leq &C\mathbb{E}\Big[|\Delta \xi|^2
+\int_t^T |\Delta f(s)|^2ds
+  \sup_{s\in[t,T]}|\Delta S_s|^2\Big|\mathcal{G}_t\Big].
\end{aligned}
\end{equation}
Step 2.  We show that
\begin{equation}\nonumber
\begin{aligned}
\mathbb{E}[ \sup_{0\leq t\leq T}|\Delta Y_t|^2]
+\mathbb{E}[\sup_{0\leq t\leq T}|\Delta K_T-\Delta K_t|^2]
\leq C\mathbb{E}\Big[|\Delta \xi|^2
+\int_t^T |\Delta f(s)|^2ds
+  \sup_{t\in[0,T]}|\Delta S_t|^2\Big].
\end{aligned}
\end{equation}
From Proposition 3.1, we have, for $i=1,2,$
\begin{equation}\nonumber
\begin{aligned}
K_T^i-K_t^i
=&\sup_{t\leq s\leq T}\Big(\mathbb{E}[\xi^i|\mathcal{G}_T]+
\mathbb{E}[\int_0^Tf^i(r,Y_r^i,Z_r^i)dr|\mathcal{G}_T]-\mathbb{E}[\int_0^{s}
f^i(r,Y_r^i,Z_r^i)dr|\mathcal{G}_{s}]\\
&\qquad\quad+\mathbb{E}[\int_0^{s}Z^i_rdW_r|\mathcal{G}_{s}]
-\mathbb{E}[\int_0^TZ^i_rdW_r|\mathcal{G}_T]
-\mathbb{E}[S_{s}^i|\mathcal{G}_{s}]
\Big)^-,
\end{aligned}
\end{equation}
which implies that
\begin{equation}\label{091604}
\begin{aligned}
|\Delta K_T-\Delta K_t|
\leq&\mathbb{E}[\big|\Delta \xi\big||\mathcal{G}_T]+\mu\mathbb{E}[\int_0^T|\Delta Y_r|+|\Delta Z_r|dr
|\mathcal{G}_T]+\mathbb{E}[\int_0^T|\Delta f(r)|dr
|\mathcal{G}_T]\\
&
+\mu\sup_{0\leq s\leq T}\mathbb{E}[\int_0^T|\Delta Y_r|+|\Delta Z_r|dr
|\mathcal{G}_s]+\sup_{0\leq s\leq T}\mathbb{E}[\int_0^T|\Delta f(r)|dr
|\mathcal{G}_s]\\
&+\sup_{0\leq s\leq T}\mathbb{E}[\sup_{0\leq s\leq T}\big|\int_0^{s}\Delta Z_rdW_r\big||\mathcal{G}_{s}]+\mathbb{E}[\big|\int_0^T\Delta Z_rdW_r\big||\mathcal{G}_T]
+\sup_{0\leq s\leq T}\mathbb{E}[\sup_{0\leq s\leq T}\big|\Delta S_{s}\big||\mathcal{G}_{s}].
\end{aligned}
\end{equation}
Thus, from \eqref{091604}, Doob's martingale inequality and Burkholder-Davis-Gundy inequality, there exists a constant $C$ only relying on $T$ and $\mu$ such that
\begin{equation}\label{091605}
\begin{aligned}
&\mathbb{E}[\sup_{0\leq t\leq T}|\Delta K_T-\Delta K_t|^2]\\
\leq&C\mathbb{E}[\big|\Delta \xi\big|^2]
+C\mathbb{E}[\int_0^T|\Delta Y_r|^2+|\Delta Z_r|^2dr]
+C\mathbb{E}[\int_0^T|\Delta f(r)|^2dr]
+C\mathbb{E}[\sup_{0\leq s\leq T}\big|\Delta S_{s}\big|^2].
\end{aligned}
\end{equation}
Then, it follows from \eqref{091603} and \eqref{091605} that
\begin{equation}\label{091606}
\begin{aligned}
\mathbb{E}[\sup_{0\leq t\leq T}|\Delta K_T-\Delta K_t|^2]
\leq C\mathbb{E}[\big|\Delta \xi\big|^2+\int_0^T|\Delta f(r)|^2dr+\sup_{0\leq t\leq T}\big|\Delta S_{t}\big|^2].
\end{aligned}
\end{equation}
On the other hand, since
\begin{equation}\nonumber
\Delta Y_t=\mathbb{E}[\Delta\xi|\mathcal{F}_t]+\mathbb{E}[\int_t^Tf^1(s,Y_s^1,Z_s^1)
-f^1(s,Y_s^2,Z_s^2)+\Delta f(s)ds|\mathcal{F}_t]
+\mathbb{E}[\Delta K_T-\Delta K_t|\mathcal{F}_t],
\end{equation}
from \eqref{091603} and \eqref{091606} we have
\begin{equation}\nonumber
\begin{aligned}
\mathbb{E}[\sup_{0\leq t\leq T}|\Delta Y_t|^2]
\leq C\mathbb{E}[\big|\Delta \xi\big|^2+\int_0^T|\Delta f(r)|^2dr+\sup_{0\leq t\leq T}\big|\Delta S_{t}\big|^2].
\end{aligned}
\end{equation}

Finally, combining Step 1 and Step 2 we get the desired result.
\end{proof}
Similar to the proof of Theorem \ref{le4.1}, we have the following result.
\begin{corollary}\label{le4.1}
Let the Assumption (H1) hold and  $(Y,Z,K)\in\mathcal{S}^2\times\mathcal{H}^2\times\mathcal{A}_{\mathbb{G}}^2$ be a solution
of conditional RBSDE \eqref{constrained-bsde}. Then there exists a constant $C$ only depending on $T$ and $\mu$ such that, for any $t\in[0,T]$,
\begin{equation}\nonumber
\mathbb{E}\Big[|Y_t|^2+\int_t^T[|Y_s|^2
+|Z_s|^2]ds+|K_T-K_t|^2\Big|\mathcal{G}_t\Big]\leq
C\mathbb{E}\Big[|\xi|^2+\int_t^T|f(s,0,0)|^2ds+\sup_{t\leq s\leq T}|S_s|^2\Big|\mathcal{G}_t\Big],\ \mathbb{P}\text{-a.s.,}
\end{equation}
and
\begin{equation}\nonumber
\mathbb{E}\Big[\sup_{0\leq t\leq T}|Y_t|^2+\int_0^T
|Z_t|^2]dt+|K_T|^2\Big]\leq
C\mathbb{E}\Big[|\xi|^2+\int_0^T|f(s,0,0)|^2ds+\sup_{0\leq t\leq T}|S_t|^2\Big].
\end{equation}
\end{corollary}
As a byproduct of Theorem \ref{le4.1}, we obtain the following uniqueness result  directly.
\begin{theorem}\label{unique}
Let the parameter $(\xi,f,S)$ satisfy the Assumption (H1). Then  conditional RBSDE
\eqref{constrained-bsde}  has at most a  solution $(Y,Z,K)$ in $\mathcal{S}^2\times\mathcal{H}^2\times\mathcal{A}_{\mathbb{G}}^2$.
\end{theorem}

\subsection{Existence of a solution}

We first focus on the particular case when $f$ do not depend on $(y,z)$, i.e., $f(s,y,z)\equiv f(s)$. In this case, we construct explicitly a solution via Snell envelope approach, i.e., an associated optimal stopping problem with partial information. For each $t\in[0,T]$, we denote by $\mathcal{T}_{t,T}$  the set of $\mathbb{G}$-adapted stopping times of values in $[t,T]$.

\begin{proposition}\label{prop3.1}
Let the parameter $(\xi,f,S)$ satisfy the Assumption (H1) and $f(s,y,z)\equiv f(s)$. Then conditional RBSDE
\eqref{constrained-bsde}  has  a unique solution $(Y,Z,K)\in\mathcal{S}^2\times\mathcal{H}^2\times\mathcal{A}_{\mathbb{G}}^2$.
\end{proposition}
\begin{proof}
For each $t\in[0,T]$ and $\tau\in\mathcal{T}_{t,T}$, we denote by $(y^\tau,z^\tau)\in\mathcal{S}^2\times\mathcal{H}^2$ the unique solution of the following BSDE
\begin{equation}\nonumber
y_s^\tau=\big[\xi I_{\{\tau=T\}}+S_\tau I_{\{\tau<T\}}\big]+\int_s^\tau f(r)dr -\int_s^\tau z_r^\tau dW_r,\ s\in[t,\tau].
\end{equation}
Then we  consider an optimal
stopping problem: for each $t\in[0,T]$,
\begin{equation}\label{050202}
\begin{aligned}
\overline{Y}_t:=\esssup_{\tau\in\mathcal{T}_{t,T}}
\mathbb{E}[y_t^\tau|\mathcal{G}_t]
=\esssup_{\tau\in\mathcal{T}_{t,T}}\mathbb{E}\Big[\big[\xi I_{\{\tau=T\}}+S_\tau I_{\{\tau<T\}}\big]+\int_t^\tau f(s)ds \Big|\mathcal{G}_t\Big].
\end{aligned}
\end{equation}
Obviously, the value process $\overline{Y}$ of the optimal stopping is $\mathbb{G}$-adapted. It follows from \eqref{050202} that
\begin{equation}\nonumber
\begin{aligned}
\overline{Y}_t+\mathbb{E}\big[\int_0^t f(s)ds |\mathcal{G}_t\big]=
\esssup_{\tau\in\mathcal{T}_{t,T}}\mathbb{E}\Big[\big[\xi I_{\{\tau=T\}}+S_\tau I_{\{\tau<T\}}\big]+\int_0^\tau f(s)ds \Big|\mathcal{G}_t\Big].
\end{aligned}
\end{equation}
Thus,
$\{\overline{Y}_t+E\big[\int_0^t f(s)ds |\mathcal{G}_t\big]\}_{t\in[0,T]}$ is the \emph{Snell envelope} of the process $\{H_t\}_{t\in[0,T]}$, where
$$H_t:=\mathbb{E}\Big[\xi I_{\{t=T\}}+S_t I_{\{t<T\}}+\int_0^t f(s)ds \Big|\mathcal{G}_t\Big],$$
that is, it is the smallest continuous $\mathbb{G}$-supermartingale that dominates the process $H$.
The continuity property of $\overline{Y}$ follows from  the fact that the process
 $H$ is continuous on $[0,T)$
and  the jump at $T$ is nonnegative.
Then it follows from the Doob-Meyer decomposition theorem that, there exists
a continuous process $K\in\mathcal{A}_\mathbb{G}^2$  and uniformly integrable $\mathbb{G}$-martingale $M$ such that
\begin{equation}\label{050203}
\overline{Y}_t+\mathbb{E}\big[\int_0^t f(s)ds |\mathcal{G}_t\big]=\mathbb{E}\big[\xi+\int_0^T f(s)ds |\mathcal{G}_T\big]+K_T-K_t-(M_T-M_t).
\end{equation}
Since $\overline{Y}$ is $\mathbb{G}$-adapted, from \eqref{050203} we get
\begin{equation}\label{050206}
\overline{Y}_t=\mathbb{E}\big[\xi+\int_t^T f(s)ds |\mathcal{G}_t\big]+\mathbb{E}[K_T-K_t|\mathcal{G}_t\big].
\end{equation}
On the other hand, it follows from \eqref{050202} and the classical optimal
stopping theory (see, for example, Proposition B.11 in \cite{KQ2012} or Theorem D.13 in \cite{KS1998}) that $\overline{Y}_t\geq \mathbb{E}[S_t|\mathcal{G}_t],\ t\in[0,T],$ and
\begin{equation}\label{050204}
\begin{aligned}
\int_0^T(\overline{Y}_t- \mathbb{E}[S_t|\mathcal{G}_t])dK_t=\int_0^T
\Big(\overline{Y}_t+\mathbb{E}[\int_0^t f(s)ds |\mathcal{G}_t]-\mathbb{E}[S_t+\int_0^t f(s)ds |\mathcal{G}_t])dK_t=0.
\end{aligned}
\end{equation}
Along with the process $K$ obtained above, the following BSDE
\begin{equation}\label{050205}
Y_t=\xi+\int_t^T f(s)ds +K_T-K_t-\int_t^TZ_sdW_s,\ t\in[0,T],
\end{equation}
has a unique solution $(Y,Z)\in\mathcal{S}^2\times\mathcal{H}^2$.
Notice that
$\overline{Y}_t=\mathbb{E}[Y_t|\mathcal{G}_t]$ because of \eqref{050206} and \eqref{050205},  combining  \eqref{050205} with \eqref{050204} we show that $(Y,Z,K)$ is a solution of conditional RBSDE \eqref{constrained-bsde}. The uniqueness follows from Theorem \ref{unique}.
\end{proof}

We now turn to the general driver case and show the existence of a solution combining Proposition \ref{prop3.2}, Proposition \ref{prop3.1} and contraction arguments.

\begin{theorem}\label{th3.1}
Suppose that the parameter $(\xi,f,S)$ satisfies the Assumption (H1). Then conditional RBSDE
\eqref{constrained-bsde}   has a unique  solution $(Y,Z,K)\in\mathcal{S}^2\times\mathcal{H}^2\times\mathcal{A}_{\mathbb{G}}^2$.
\end{theorem}
\begin{proof}
We only need to prove the existence of a solution since the uniqueness has been obtained in
Theorem  \ref{unique}.
For any given $(U,V)\in\mathcal{S}^2\times\mathcal{H}^2$, it follows from Proposition
\ref{prop3.1} that the following equation
\begin{equation}\label{050211}
\left\{
\begin{aligned}
&Y_t=\xi+\int_t^T f(s,U_s,V_s)ds +K_T-K_t-\int_t^TZ_sdW_s,\ t\in[0,T],\\
&\mathbb{E}[Y_t-S_t|\mathcal{G}_t]\geq 0,\ \forall\ t\in[0,T],\  \int_0^T\mathbb{E}[Y_t-S_t|\mathcal{G}_t]dK_t=0,
\end{aligned}
\right.
\end{equation}
exists a unique solution $(Y,Z,K)\in\mathcal{S}^2\times\mathcal{H}^2\times\mathcal{A}_{\mathbb{G}}^2$.
Moreover,  using Proposition \ref{prop3.2} it holds
\begin{equation}\label{082902}
\begin{aligned}
K_T-K_t
=&\sup_{t\leq s\leq T}\Big(\mathbb{E}[\xi|\mathcal{G}_T]+
\mathbb{E}[\int_0^Tf(r,U_r,V_r)dr|\mathcal{G}_T]-\mathbb{E}[\int_0^{s}
f(r,U_r,V_r)dr|\mathcal{G}_{s}]\\
&\qquad\quad+\mathbb{E}[\int_0^{s}Z_rdW_r|\mathcal{G}_{s}]
-\mathbb{E}[\int_0^TZ_rdW_r|\mathcal{G}_T]
-\mathbb{E}[S_{s}|\mathcal{G}_{s}]
\Big)^-.
\end{aligned}
\end{equation}
Thus, using \eqref{082902} and \eqref{050211} we may define a mapping from Banach space $\mathcal{S}^2\times\mathcal{H}^2$ to itself  as
$$\Phi:(U,V)\rightarrow (Y,Z),$$
and only need to show that it is a  contraction mapping. For $(U^1,V^1),(U^2,V^2)\in\mathcal{S}^2\times\mathcal{H}^2$, we denote
$$(Y^1,Z^1)=\Phi(U^1,V^1),\ (Y^2,Z^2)=\Phi(U^2,V^2);\
\Delta L=L^1-L^2,\ L=Y,Z,U,V,K.$$
  Classical arguments suggest that, for any $\beta>0$, we have
\begin{equation}\label{050213}
\begin{aligned}
&\mathbb{E}[e^{\beta t}|\Delta Y_t|^2|\mathcal{G}_t]
+\mathbb{E}[\int_t^T\beta e^{\beta s}|\Delta Y_s|^2 ds
+\int_t^T e^{\beta s}|\Delta Z_s|^2ds|\mathcal{G}_t]\\
=&\mathbb{E}[\int_t^T2 e^{\beta s}\Delta Y_s\big(f(s,U_s^1,V_s^1)-f(s,U_s^2,V_s^2)\big)ds|\mathcal{G}_t]
+\mathbb{E}[\int_t^T2 e^{\beta s}\Delta Y_sd\Delta K_s|\mathcal{G}_t]\\
\leq &\mathbb{E}[\int_t^T2 e^{\beta s}\Delta Y_s\big(f(s,U_s^1,V_s^1)-f(s,U_s^2,V_s^2)\big)ds|\mathcal{G}_t]\quad (\text{from \eqref{110401} with}\ \Delta S\equiv 0)\\
\leq &4\mu^2\mathbb{E}[\int_t^T e^{\beta s}|\Delta Y_s|^2ds|\mathcal{G}_t]+\frac{1}{2} \mathbb{E}[\int_t^T e^{\beta s}\big(|\Delta U_s|^2+|\Delta V_s|^2\big)ds|\mathcal{G}_t].
\end{aligned}
\end{equation}
By choosing $\beta=4\mu^2+1$, we have
\begin{equation}\label{050216}
\begin{aligned}
\mathbb{E}[e^{\beta t}|\Delta Y_t|^2|\mathcal{G}_t]
+\mathbb{E}[\int_t^T e^{\beta s}|\Delta Y_s|^2 ds
+\int_t^T e^{\beta s}|\Delta Z_s|^2ds|\mathcal{G}_t]
\leq \frac{1}{2} \mathbb{E}[\int_t^T e^{\beta s}\big(|\Delta U_s|^2+|\Delta V_s|^2\big)ds|\mathcal{G}_t],
\end{aligned}
\end{equation}
which implies that $\Phi$ is a strict contraction mapping on $\mathcal{H}^2\times\mathcal{H}^2$ with the norm
$
\|(Y,Z)\|_{\beta}^2=\mathbb{E}\int_0^Te^{\beta t}(|Y_t|^2+|Z_t|^2)dt.
$
 {On the other hand, from \eqref{082902} we have
\begin{equation}\label{082903}
\begin{aligned}
|\Delta K_T-\Delta K_t|
\leq& \mu\mathbb{E}\Big[\int_0^T(|\Delta U_r|+|\Delta V_r|)dr\Big|\mathcal{G}_T\Big]
+\mu\sup_{0\leq s\leq T}\mathbb{E}\Big[\int_0^T(|\Delta U_r|+|\Delta V_r|)dr\Big|\mathcal{G}_s\Big]\\
&+\sup_{0\leq s\leq T}\mathbb{E}\Big[\sup_{0\leq s\leq T}|\int_0^s\Delta Z_rdW_r|\Big|\mathcal{G}_s\Big]
+\mathbb{E}\Big[|\int_0^T\Delta Z_rdW_r|\Big|\mathcal{G}_T\Big].
\end{aligned}
\end{equation}
Since
$\Delta Y_t=\mathbb{E}[\int_t^Tf(s,U_s^1,V_s^1)-f(s,U_s^2,V_s^2)ds|\mathcal{F}_t]
+\mathbb{E}[\Delta K_T-\Delta K_t|\mathcal{F}_t],$
from \eqref{082903} we have
\begin{equation}\label{090601}
\begin{aligned}
&|\Delta Y_t|
\leq \mu\mathbb{E}\Big[\int_0^T(|\Delta U_s|+|\Delta V_s|)ds\Big|\mathcal{F}_t\Big]
+\mu\mathbb{E}\Big[\mathbb{E}[\int_0^T(|\Delta U_s|+|\Delta V_s|)ds|\mathcal{G}_T]\Big|\mathcal{F}_t\Big]\\
&+\mu\mathbb{E}\Big[\sup_{0\leq s\leq T}\mathbb{E}[\int_0^T(|\Delta U_r|+|\Delta V_r|)dr|\mathcal{G}_s]\Big|\mathcal{F}_t\Big]
+\mathbb{E}\Big[\sup_{0\leq s\leq T}\mathbb{E}[\sup_{0\leq s\leq T}|\int_0^s\Delta Z_rdW_r|\big|\mathcal{G}_s]\Big|\mathcal{F}_t\Big]\\
&+\mathbb{E}\Big[\mathbb{E}[|\int_0^T\Delta Z_sdW_s||\mathcal{G}_T]\Big|\mathcal{F}_t\Big].
\end{aligned}
\end{equation}
Then it follows from Doob's martingale inequality, \eqref{090601} and \eqref{050216} that
\begin{equation}\nonumber
\mathbb{E}[\sup_{0\leq t\leq T}|\Delta Y_t|^2]\leq C\mathbb{E}[\int_0^T|\Delta U_s|^2+|\Delta V_s|^2ds].
\end{equation}
As a result, $\Phi$ is continuous from $\mathcal{S}^2\times\mathcal{H}^2$ to itself. Combining with \eqref{050216},  $\Phi$ has a unique fixed point $(Y,Z)\in\mathcal{S}^2\times\mathcal{H}^2$. The existence of $K$ follows  directly from \eqref{050211} and \eqref{082902}.
}
\end{proof}

\begin{remark}
We obtain the well-posedness of the solution of conditional RBSDE \eqref{constrained-bsde} when the reflection condition is linear in $y$. It seems more
interesting to consider such equation with nonlinear reflection condition similar to the study of BSDEs with general mean reflection in \cite{BEH2018} (see, Section 4 therein). In fact, we can extend Theorem \ref{th3.1} to the general nonlinear conditional reflection situation
$$\ell(t,\mathbb{E}[Y_t|\mathcal{G}_t])\geq 0,\ 0\leq t\leq T,$$
with some increasing (in $y$) continuous function $\ell$. The proof is similar to Theorem 9 in \cite{BEH2018} by introducing an operator $L_t$ defined as
$$L_t:\mathcal{L}^2(\mathcal{G}_T)\mapsto
\mathcal{L}^2(\mathcal{G}_T;\mathbb{R}^+),\
L_t(X)=\essinf\big\{\eta:\eta\geq 0, a.s., \mathcal{G}_T\text{-measurable}, \ell(t,X+\eta)\geq 0\big\}.$$
For the study of BSDEs with another general nonlinear conditional reflection   $\mathbb{E}[\ell(t,Y_t)|\mathcal{G}_t]\geq 0,\ 0\leq t\leq T,$
we will leave it for the future  research.
\end{remark}

\section{The connection between optimal stopping problems and linear conditional RBSDEs}

In this section, we study the connection between conditional RBSDE \eqref{constrained-bsde} and the related  optimal stopping problems when the driver $f(t,y,z)$ is linear in $(y,z)$. First of all, from the proof of Proposition \ref{prop3.1}, we get the following connection  when the driver  $f$ does not depend on $(y,z)$, i.e., $f(s,y,z)\equiv f(s)$.

\begin{corollary}\label{co4.2}
Let $(Y,Z,K)$ be the solution of the following conditional RBSDE
\begin{equation}\nonumber
\left\{
\begin{aligned}
&Y_t=\xi+\int_t^Tf(s)ds +K_T-K_t-\int_t^TZ_s dW_s,\ t\in[0,T],\\
&\mathbb{E}[Y_t-S_t|\mathcal{G}_t]\geq 0,\ \forall\ t\in[0,T],\  \int_0^T\mathbb{E}[Y_t-S_t|\mathcal{G}_t]dK_t=0.
\end{aligned}
\right.
\end{equation}
Then, we have, for all $t\in[0,T]$,
\begin{equation}\nonumber
\mathbb{E}[Y_t|\mathcal{G}_t]=\esssup_{\tau\in\mathcal{T}_{t,T}}
\mathbb{E}[y_t^\tau|\mathcal{G}_t],
\end{equation}
where, for each $\tau\in\mathcal{T}_{t,T}$,  $(y^\tau,z^{\tau})$ is the unique solution  of the following BSDE
\begin{equation}\nonumber
y_s^\tau=\big[\xi I_{\{\tau=T\}}+S_\tau I_{\{\tau<T\}}\big]+\int_s^\tau f(r)dr -\int_s^\tau z_r^{\tau}dW_r,\ s\in[t,\tau].
\end{equation}
Moreover, the  optimal  stopping $\tau^*\in\mathcal{T}_{t,T}$ is given by
$
\tau^*_t=\inf\big\{s\in[t,T]: \mathbb{E}[Y_s-S_s|\mathcal{G}_s]=0\big\}\wedge T.
$
\end{corollary}
In particular, if  partial information $\mathbb{G}$ is chosen to be deterministic scenario, i.e., $\mathcal{G}_t=\mathcal{H}$, $0\leq t\leq T$, in Corollary \ref{co4.2}, then we get the following  connection between BSDEs with mean reflection and deterministic stopping time problems. In this case, the set of $\mathbb{G}$-adapted stopping times $\mathcal{T}_{t,T}=[t,T]$.
\begin{example}\label{eg4.1}
Let $(Y,Z,K)$ be the solution of the following BSDE with mean reflection
\begin{equation}\nonumber
\left\{
\begin{aligned}
&Y_t=\xi+\int_t^Tf(s)ds +K_T-K_t-\int_t^TZ_s dW_s,\ t\in[0,T],\\
&\mathbb{E}[Y_t-S_t]\geq 0,\ \forall\ t\in[0,T],\ \int_0^T\mathbb{E}[Y_t-S_t]dK_t=0.
\end{aligned}
\right.
\end{equation}
 Then we have, for all $t\in[0,T]$,
\begin{equation}\nonumber
\mathbb{E}[Y_t]=\sup_{\tau\in[t,T]}
\mathbb{E}[y_t^\tau],
\end{equation}
where, for each $\tau\in[t,T]$,  $(y^\tau,z^{\tau})$ is the unique solution  of the following BSDE
\begin{equation}\nonumber
y_s^\tau=\big[\xi I_{\{\tau=T\}}+S_\tau I_{\{\tau<T\}}\big]+\int_s^\tau f(r)dr -\int_s^\tau z_r^{\tau}dW_r,\ s\in[t,\tau],
\end{equation}
and an  optimal  time $\tau^*\in[t,T]$ is given by
$
\tau^*_t=\inf\big\{s\in[t,T]: \mathbb{E}[Y_s]=\mathbb{E}[S_s]\big\}\wedge T.
$
\end{example}
We now generalize Corollary \ref{co4.2} to the linear driver case, in which we specify a class of partial information $\mathbb{G}$.
For simplicity, the underlying Brownian motion is chosen to be two-dimensional, i.e., $W=(W^1,W^2)$ and recall that the filtration $\mathbb{F}$ is generated by $W$. Let
$$\mathbb{U}=
\left(
\begin{array}{cc}
\lambda_1&\lambda_2\\
\lambda_3&\lambda_4
\end{array}
\right),$$
be a constant orthogonal matrix (i.e., $\mathbb{U}\mathbb{U}^T$ is the identity matrix). Then the process $\widetilde{W}=(\widetilde{W}^1,\widetilde{W}^2)$ defined as
\begin{equation}\nonumber
(\widetilde{W}^1,\widetilde{W}^2)^T=\mathbb{U}\cdot(W^1,W^2)^T=(\lambda_1 W^1+\lambda_2 W^2, \lambda_3W^1+\lambda_4W^2)^T,
\end{equation}
is also a Brownian motion. It is easy to check that the filtration generated by $\widetilde{W}$ is still $\mathbb{F}$.
 Let $\mathbb{G}$ be the
sub-filtration generated by $\widetilde{W}^1$ (i.e., $\lambda_1 W^1+\lambda_2 W^2$).
\begin{remark}\label{re4.1}
When the Brownian motion $W$ is $d$-dimensional, one can similarly construct a new Brownian motion $\widetilde{W}$ through a $d\times d$ constant orthogonal matrix. Then the sub-filtration $\mathbb{G}$ is generated by some components of this new Brownian motion $\widetilde{W}$. We consider two dimensional situation only to simplify the notations.
On the other hand, it is easy to check that the number of such orthogonal matrix is infinity, which means that our results can be applied to a large class of partial information problems.
\end{remark}
Suppose that the driver $f$ has the following linear form
\begin{equation}\nonumber
f(s,y,z^1,z^2)=a_sy+\lambda_1b_sz^1+\lambda_2b_sz^2+c_s
=a_sy+b_s(\lambda_1,\lambda_2)\cdot (z^1,z^2)^T+c_s,
\end{equation}
where both $a$ and $b$ are $\mathbb{G}$-adapted and bounded processes, the process $c$ is $\mathbb{F}$-adapted and belongs to $\mathcal{H}^2$.
\begin{theorem}\label{le3.3}
Let $(Y,Z^1,Z^2,K)$ be the unique solution of the following conditional RBSDE
\begin{equation}\label{eq6}
\left\{
\begin{aligned}
&Y_t=\xi+\int_t^T[a_sY_s+b_s(\lambda_1,\lambda_2)\cdot (Z^1_s,Z^2_s)^T+c_s]ds +K_T-K_t-\int_t^T(Z_s^1,Z_s^2)\cdot d(W^1,W^2)_s^T,\ t\in[0,T],\\
&\mathbb{E}[Y_t-S_t|\mathcal{G}_t]\geq 0,\ \forall\ t\in[0,T],\  \int_0^T\mathbb{E}[Y_t-S_t|\mathcal{G}_t]dK_t=0.
\end{aligned}
\right.
\end{equation}
Then, we have, for all $t\in[0,T]$,
\begin{equation}\label{3}
\mathbb{E}[Y_t|\mathcal{G}_t]=\esssup_{\tau\in\mathcal{T}_{t,T}}
\mathbb{E}[y_t^\tau|\mathcal{G}_t],
\end{equation}
where, for each $\tau\in\mathcal{T}_{t,T}$,  $(y^\tau,z^{1,\tau},z^{2,\tau})$ is the unique solution  of the following BSDE: for $s\in[t,\tau]$,
\begin{equation}\label{eq7}
y_s^\tau=\big[\xi I_{\{\tau=T\}}+S_\tau I_{\{\tau<T\}}\big]+\int_s^\tau [a_r y_r^\tau+b_r(\lambda_1,\lambda_2)\cdot (z^{1,\tau}_r, z^{2,\tau}_r)^T +c_r]dr -\int_s^\tau (z_r^{1,\tau},z_r^{2,\tau})\cdot d(W^1,W^2)^T_t.
\end{equation}
Moreover, an optimal  stopping $\tau^*\in\mathcal{T}_{t,T}$ is given by
\begin{equation}\label{5}
\tau^*_t=\inf\big\{s\in[t,T]: \mathbb{E}[Y_s-S_s|\mathcal{G}_s]=0\big\}\wedge T.
\end{equation}
\end{theorem}

\begin{proof}
 Let $t\in[0,T]$ be arbitrarily fixed.
In order to show \eqref{3} and \eqref{5}, we only need to show that\\
(i) For all $\tau\in\mathcal{T}_{t,T}$, $\mathbb{E}[Y_t|\mathcal{G}_t]\geq
\mathbb{E}[y_t^\tau|\mathcal{G}_t]$;\\
(ii) With $\tau^*_t$ given in \eqref{5}, it holds $\mathbb{E}[Y_t|\mathcal{G}_t]=
\mathbb{E}[y_t^{\tau^*_t}|\mathcal{G}_t]$.\\
\noindent (i) For each $\tau\in\mathcal{T}_{t,T}$, from \eqref{eq6} and \eqref{eq7}, $(\Delta Y,\Delta Z^1,\Delta Z^2)$
$:=(Y-y^\tau,Z^1-z^{1,\tau},Z^2-z^{2,\tau})$ satisfies
\begin{equation}\nonumber
\Delta Y_t=\Delta Y_{\tau}+\int_t^{\tau} [a_s\Delta Y_s+b_s(\lambda_1,\lambda_2)\cdot (\Delta Z^{1}_s, \Delta Z^{2}_s)^T]ds+K_\tau-K_t -\int_t^{\tau}(\Delta Z_s^1,\Delta Z_s^2)\cdot d(W^1,W^2)_s^T,
\end{equation}
with the terminal condition $\Delta Y_\tau=Y_\tau-\big[\xi I_{\{\tau=T\}}+S_\tau I_{\{\tau<T\}}\big].$
Let $\Gamma$ be the unique solution of the following SDE
\begin{equation}\label{eq9}
\begin{aligned}
d\Gamma_s=a_s\Gamma_sds+b_s\Gamma_s(\lambda_1,\lambda_2)d(W^1,W^2)^T_s
=a_s\Gamma_sds+b_s\Gamma_sd\widetilde{W}^1_s,\ s\in[t,T];\
\Gamma_t=1.
\end{aligned}
\end{equation}
Then we have
$
\Gamma_s
=\exp\{\int_t^s(a_r-\frac12 b_r^2)dr+\int_t^sb_rd\widetilde{W}_r^1\}\in\mathcal{G}_s,\ s\in[t,T].
$ 
Using It\^o's formula to $\Gamma_s \Delta Y_s$, we have
\begin{equation}\label{eq11}
\Delta Y_t=\mathbb{E}[\Gamma_\tau\Delta Y_{\tau}+\int_t^\tau \Gamma_sdK_s|\mathcal{F}_t]\geq \mathbb{E}[\Gamma_\tau\Delta Y_{\tau}|\mathcal{F}_t],
\end{equation}
which implies that
\begin{equation}\label{eq12}
\mathbb{E}[\Delta Y_t|\mathcal{G}_t]\geq \mathbb{E}\big[\Gamma_\tau\Delta Y_{\tau}|\mathcal{G}_t\big]
=\mathbb{E}\big[\mathbb{E}[\Gamma_\tau\Delta Y_{\tau}|\mathcal{G}_\tau]|\mathcal{G}_t\big]=\mathbb{E}\big[\Gamma_\tau\mathbb{E}[\Delta Y_{\tau}|\mathcal{G}_\tau]|\mathcal{G}_t\big].
\end{equation}
Notice that $\mathbb{E}[Y_s-S_s|\mathcal{G}_s]\geq 0$, $s\in[0,T]$,  we have $\mathbb{E}[Y_s-\xi I_{\{s=T\}}-S_sI_{\{s<T\}}|\mathcal{G}_s]\geq 0$,  from which it holds
 $${\mathbb{E}[\Delta Y_{\tau}|\mathcal{G}_\tau]
=\mathbb{E}[Y_{\tau}-\xi I_{\{\tau=T\}}-S_\tau I_{\{\tau<T\}}|\mathcal{G}_\tau]\geq 0}.$$
Then it follows from \eqref{eq12} that $\mathbb{E}[\Delta Y_t|\mathcal{G}_t]\geq0$, i.e.,
$\mathbb{E}[Y_t|\mathcal{G}_t]\geq
\mathbb{E}[y_t^\tau|\mathcal{G}_t]$. \\
\noindent (ii)
From \eqref{5} we can check
\begin{equation}\label{eq15}
\mathbb{E}[\Delta Y_{{\tau_t^*}}|\mathcal{G}_{\tau_t^*}]
=
\mathbb{E}[\big(Y_{{\tau_t^*}}-S_{\tau_t^*}\big)\cdot I_{\{{\tau_t^*}<T\}}
+\big(Y_{{\tau_t^*}}-\xi\big)\cdot I_{\{{\tau_t^*}=T\}}|\mathcal{G}_{\tau_t^*}]= 0.
\end{equation}
Combining \eqref{eq15}, noting that inequalities  \eqref{eq11} and \eqref{eq12} with $\tau=\tau_t^*$ turn to be equalities since $K_{\tau_t^*}=K_t$, we get $\mathbb{E}[Y_t|\mathcal{G}_t]= \mathbb{E}[y_t^{\tau^*_t}|\mathcal{G}_t]$.
\end{proof}

\begin{remark}
From Corollary \ref{co4.2} and Theorem \ref{le3.3}, we conclude that the link between conditional RBSDEs and optimal stopping problems can be obtained in two special cases:\\
(1) When $f(s,y,z)\equiv f(s)$,  partial information $\mathbb{G}$ needs no  requirement  except the basic assumption;\\
(2) When $f$ is linear in $(y,z)$, partial information $\mathbb{G}$  may need some specific structure as given above. \\
The study of the link between conditional RBSDEs and optimal stopping problems with general nonlinear driver $f$ and partial information $\mathbb{G}$ is left  for the further research.
\end{remark}
With the help of Theorem \ref{le3.3}, we can show that comparison theorem holds for linear conditional RBSDE \eqref{eq6}.
\begin{corollary}\label{Co4.1}(Comparison Theorem)
Suppose that $(\xi_i,S^i),$ $i=1,2,$ satisfy the Assumption (H1).
Let $(Y^i,Z^i,K^i)$, $i=1,2$, be the unique solution of conditional RBSDE
\begin{equation}\nonumber
\left\{
\begin{aligned}
&Y_t^i=\xi_i+\int_t^T[a_sY_s^i+b_s(\lambda_1,\lambda_2)Z^i_s+c_s^i]ds +K_T^i-K_t^i-\int_t^TZ_s^i dW_s,\ t\in[0,T],\\
&\mathbb{E}[Y_t^i-S_t^i|\mathcal{G}_t]\geq 0,\ \forall\ t\in[0,T];\  \int_0^T\mathbb{E}[Y_t^i-S_t^i|\mathcal{G}_t]dK_t^i=0.
\end{aligned}
\right.
\end{equation}
If the following conditions hold:\\
\textup{(1)} For the terminal conditions $\xi_1$ and $\xi_2$, $\mathbb{E}[\xi_1|\mathcal{G}_T]\geq \mathbb{E}[\xi_2|\mathcal{G}_T]$,\\
\textup{(2)} For the processes $c^1$ and $c^2$, $\mathbb{E}[c^1_t|\mathcal{G}_t]\geq \mathbb{E}[c^2_t|\mathcal{G}_t]$, for $t\in[0,T]$,\\
\textup{(3)} For the barriers $S^1$ and $S^2$, $\mathbb{E}[S^1_t|\mathcal{G}_t]\geq \mathbb{E}[S^2_t|\mathcal{G}_t]$, for $t\in[0,T]$,\\
then for each $t\in[0,T]$,  we get
$$\mathbb{E}[Y_t^1|\mathcal{G}_t]\geq \mathbb{E}[Y_t^2|\mathcal{G}_t],\ \mathbb{P}\text{-a.s.}$$
\end{corollary}

\begin{proof}
From Theorem \ref{le3.3}, we get, for $t\in[0,T]$, $i=1,2$,
\begin{equation}\label{070803}
\mathbb{E}[Y_t^i|\mathcal{G}_t]
=\esssup_{\tau\in\mathcal{T}_{t,T}}\mathbb{E}[y_t^i|\mathcal{G}_t],
\end{equation}
where $(y^i,z^i)$ is the unique solution of the following BSDE
\begin{equation}\nonumber
y_s^i=[\xi_iI_{\{\tau=T\}}+S_\tau^i I_{\{\tau<T\}}]+\int_s^\tau [a_r y_r^i+b_r(\lambda_1,\lambda_2) z^{i}_r +c_r^i]dr -\int_s^\tau z_r^{i} dW_r,\ s\in[t,\tau].
\end{equation}
Similar to the proof  of Theorem \ref{le3.3}, we get
\begin{equation}\label{070801}
\mathbb{E}[\Delta y_t|\mathcal{G}_t]=\mathbb{E}[\Gamma_\tau\Delta y_{\tau}+\int_t^\tau \Gamma_s\Delta c_sds|\mathcal{G}_t],
\end{equation}
where $\Gamma$ is the solution of SDE \eqref{eq9}, $\Delta y_t=y_t^1-y_t^2$, $\Delta y_\tau=(\xi_1-\xi_2)I_{\{\tau=T\}}+(S_\tau^1-S_\tau^2)I_{\{\tau<T\}}$ and $\Delta c_s=c_s^1-c_s^2$.
Since
\begin{equation}\nonumber
\mathbb{E}[\Gamma_\tau\Delta y_{\tau}|\mathcal{G}_t]
=\mathbb{E}[\Gamma_\tau\mathbb{E}[\Delta y_{\tau}|\mathcal{G}_\tau]|\mathcal{G}_t]\geq 0,\
\mathbb{E}[\int_t^\tau \Gamma_s\Delta c_sds|\mathcal{G}_t]
=\mathbb{E}[\int_t^\tau \Gamma_s\mathbb{E}[\Delta c_s|\mathcal{G}_s]ds|\mathcal{G}_t]\geq 0,
\end{equation}
it follows from \eqref{070801} that
$\mathbb{E}[y_t^1|\mathcal{G}_t]\geq \mathbb{E}[y_t^2|\mathcal{G}_t],\ \text{for\ each}\ \tau\in\mathcal{T}_{t,T},$
from which we can conclude the desired result by using \eqref{070803}.
\end{proof}
\begin{remark}
As the Example 3.3 in \cite{HHLLW} shows that we can not expect to compare $Y^1$ and $Y^2$ pointwisely in Corollary \ref{Co4.1}. It seems reasonable to compare these two terms under the conditional expectation with respect to the partial information $\mathbb{G}$.
\end{remark}

\section{Backward recursive reflected control problems with partial information}
In this section, inspired by \emph{Example 1.2} in Introduction, we  consider backward recursive reflected control problems with partial information (BRR problems, for short), where
the payoff is given by  controlled conditional RBSDEs.
Throughout this section, we adopt the partial information $\mathbb{G}$ introduced in Section 3.
For simplicity of notation, we choose the orthogonal matrix $\mathbb{U}$ to be an identity matrix (i.e., $\lambda_1=\lambda_4=1,$ $\lambda_2=\lambda_3=0$) and thus the partial information $\mathbb{G}$ is generated by the first component $W^1$ of the Brownian motion $W=(W^1,W^2)$.
Let $V$ be a nonempty compact subset of $\mathbb{R}^k$.
An \emph{admissible} control $v:[0,T]\times\Omega\rightarrow V$ is an $\mathbb{G}$-adapted process {such that}   $\mathbb{E}\Big[\int_0^T|v_t|^2dt\Big]<\infty,$ and we denote by $\mathcal{V}$ the set of all admissible controls. Herein,  $\mathbb{G}$ represents the information available to the controller, which is usually incomplete in most situations.

We consider two types of weak formulations of BRR problems: linear case and convex case, for both the state equation  is described by
the  following SDE
\begin{equation}\label{031801}
\begin{aligned}
X_t=x_0+\int_0^t \sigma(s,{X_s})dW_s,\ t\in[0,T],\ x_0\in\mathbb{R}^2,\\
\end{aligned}
\end{equation}
where the coefficient $\sigma:[0,T]\times\mathbb{R}^2\rightarrow\mathbb{R}^{2\times 2}$ is Lipschitz in $x$ and  $\sigma(t,0)$ is uniformly bounded with respect to $t\in[0,T]$.
It is well known that SDE \eqref{031801}
 has a unique solution $X\in\mathcal{S}^2$.
\begin{remark}
If we choose
$$
\sigma(t,(x_1,x_2)^T)=
\left(
\begin{array}{cc}
\sigma_1 x_1 & 0\\
0 & \sigma_2 x_2
\end{array}
\right),
$$
with $\sigma_1$ and $\sigma_2$ are two given constants, then SDE \eqref{031801} can be applied to model the price of two stocks with zero return rate (see equation \eqref{101601}), namely,  $X_t=(X^1_t,X^2_t)$ stands for the price of the first and second stocks at time $t$. In this situation, $\mathbb{G}$ represents the price information of the first stock, which is assumed to be the only one announced to public market.

On the other hand, the structure of this partial information $\mathbb{G}$ can be linked to the large-population problems (see, e.g., \cite{BFH2021,BFY2013,G2016,HWW2016}) with representing the information of {common noise}.
\end{remark}
Let the function
$b:[0,T]\times\mathbb{R}^2\times V\rightarrow\mathbb{R}
$ be uniformly bounded and continuous with respect to $v$.
For each given $v\in\mathcal{V}$, we define a probability measure $\mathbb{P}^v$ on $(\Omega,\mathcal{G}_T)$, which is equivalent to $\mathbb{P}$ and whose density function is given by
\begin{equation}\label{031802}
\frac{d\mathbb{P}^v}{d\mathbb{P}}\Big|_{\mathcal{G}_T}
=\exp\Big\{\int_0^T
b(t,{\mathbb{E}[X_t|\mathcal{G}_t]},v_t)
dW_t^1
-\frac12\int_0^T|
b(t,{\mathbb{E}[X_t|\mathcal{G}_t]},v_t)|^2dt\Big\}.
\end{equation}
Herein, we assume that the controller will use the probability measure $\mathbb{P}^v$ instead of $\mathbb{P}$ to measure the performance of the related payoffs.
It seems natural that the probability measure $\mathbb{P}^v$ chosen by the controller should rely on both the available information $\mathbb{G}$ and
 the conditional unbiased estimate of the state $\mathbb{E}[X_t|\mathcal{G}_t]$ based on this information.
\begin{remark}
When considering the weak formulations of optimal control and game problems with full information,
it is common to introduce the probability measure $\mathbb{P}^v$ similar to \eqref{031802}, such as in \cite{EH2003} and \cite{HL1995}. However, such structure \eqref{031802} with partial information  is still totally new.
\end{remark}

Thanks to Girsanov Theorem, the process
\begin{equation}\label{032207}
dW^v_t:=
-\left(
\begin{array}{c}
b(t,{\mathbb{E}[X_t|\mathcal{G}_t]},
v_t)\\
0
\end{array}
\right)
dt+dW_t,\ t\in[0,T],
 \end{equation}
is a  Brownian motion under the probability measure $\mathbb{P}^v$.
Next, we introduce the payoffs of linear and convex BRR problems in Subsection 4.1 and Subsection 4.2, respectively.

\subsection{Weak formulation of linear BRR problems}

We first introduce the following linear conditional RBSDE
\begin{equation}\label{070805}
\left\{
\begin{aligned}
&Y_t^{v}=\Phi(X_T)+\int_t^T[\alpha_sY_s^{v}
+\beta_s Z^{1,v}_s+g(s,{\mathbb{E}[X_s|\mathcal{G}_s]},v_s)]ds +K_T^{v}-K_t^{v}-\int_t^T(Z_s^{1,v},Z_s^{2,v}){dW_s^v},\ t\in[0,T],\\
&\mathbb{E}[Y_t^{v}- h(t,X_t)|\mathcal{G}_t]\geq 0,\ \forall\ t\in[0,T],\ \text{a.s.},\
\int_0^T\mathbb{E}[Y_t^{v}- h(t,X_t)|\mathcal{G}_t]dK_t^{v}=0,
\end{aligned}
\right.
\end{equation}
where the mappings
$$\alpha:\Omega\times[0,T]\rightarrow\mathbb{R},\ \beta:\Omega\times[0,T]\rightarrow\mathbb{R},\ g:\Omega\times[0,T]\times\mathbb{R}^2\times V\rightarrow\mathbb{R},\ h:\Omega\times[0,T]\times\mathbb{R}^2\rightarrow\mathbb{R},\ \Phi:\Omega\times\mathbb{R}^2\rightarrow\mathbb{R},$$
are measurable and satisfy the following condition
 \begin{equation}\nonumber
{\bf (H2)}
\left\{
\begin{array}{l}
\text{(i)}\ \text{The processes}\ \alpha\ \text{and}\ \beta\ \text{are}\ \mathbb{G}\text{-adapted and uniformly bounded};\\
\text{(ii)}\ \text{For\ each}\ (x,v),\  g(\cdot,x,v)\  \text{is}\ \mathbb{G}\text{-adapted};\ g\ \text{is\ continuous\ in}\ v\ \text{and\ satisfies}\\
\qquad\qquad\qquad |g(t,x,v)|\leq C(1+|x|),\ \mathbb{P}\text{-a.s.},\ \text{for\ all}\ t\in[0,T];\\
\text{(iii)}\ \text{For\ each}\ x, h(\cdot,x)\ \text{is}\ \mathbb{F}\text{-adapted};  h\ \text{is\ continuous in}\ (t,x)\ \text{satisfying}\\
\qquad\qquad\qquad {|h(t,x)|\leq C(1+|x|)},\ \mathbb{P}\text{-a.s.},\ \text{for\ all}\ t\in[0,T]; \\ 
\text{(iv)}\ \text{For\ each}\ x,\ \Phi(x) \text{\ is}\ \mathcal{F}_T\text{-measurable}; \Phi\ \text{is\ continuous\ in}\ x\  \text{and\ satisfies}\\
\qquad\qquad\qquad  {|\Phi(x)|\leq C(1+|x|)},\  \mathbb{E}[h(T,x)|\mathcal{G}_T]\leq \mathbb{E}[\Phi(x)|\mathcal{G}_T],\ x\in\mathbb{R}^2.
 \end{array}
 \right.
 \end{equation}
For each given $v\in\mathcal{V}$, we can check from Theorem \ref{th3.1} and \eqref{032207} that conditional RBSDE \eqref{070805}  has a unique solution $(Y^v,Z^v,K^v)\in\mathcal{S}^2\times\mathcal{H}^2
\times\mathcal{A}_{\mathbb{G}}^2$ under the Assumption (H2).
%
%
The payoff of linear BRR problem with admissible control $v$ is defined as $Y^v_0$ and
the aim  is to maximize this recursive payoff  over all admissible controls $v\in\mathcal{V}$, i.e.,
\begin{equation}\label{070806}
{ \text{\bf(Linear BRR)}}\ \sup_{v\in\mathcal{V}}Y^v_0.
\end{equation}
From \eqref{070805}, we see that $Y_0^v\geq h(0,x_0)$, for all $v\in\mathcal{V}$, which implies that linear BRR problem \eqref{070806} is a type of optimization problems with an inequality-type constraint.

From Theorem \ref{le3.3}, it holds that  $Y_0^v=\sup_{\tau\in\mathcal{T}_{0,T}}\overline{y}_0^{\tau,v}$.  As a result, we only need to consider the following  mixed control problem, which is equivalent to linear BRR problem \eqref{070806},
$$ \sup_{\tau\in\mathcal{T}_{0,T}}\sup_{v\in\mathcal{V}}{\overline{y}_0^{\tau,v}},$$
where, for each $(v,\tau)\in\mathcal{V}\times\mathcal{T}_{0,T}$, $(\overline{y}^{\tau,v},\overline{z}^{\tau,v})\in\mathcal{S}^2\times\mathcal{H}^2$ is the unique solution of BSDE
 \begin{equation}\nonumber
\left\{
\begin{aligned}
-d\overline{y}_t^{\tau,v}=& \Big[\alpha_t\overline{y}_t^{\tau,v}
+\beta_t\overline{z}_t^{1,\tau,v}
+\overline{z}_t^{1,\tau,v}
\cdot b(t,{\mathbb{E}[X_t|\mathcal{G}_t]},v_t)
+g(t,{\mathbb{E}[X_t|\mathcal{G}_t]},v_t)\Big]dt - (\overline{z}_t^{1,\tau,v},\overline{z}_t^{2,\tau,v})dW_t,\ t\in[0,\tau],\\
\overline{y}_\tau^{\tau,v}=&\Big[\Phi(X_T)I_{\{\tau=T\}}
+{h}(\tau,X_\tau)I_{\{\tau<T\}}\Big],
\end{aligned}
\right.
\end{equation}
which can be rewritten as the following filtered BSDE
 \begin{equation}\nonumber
\left\{
\begin{aligned}
-d\mathbb{E}[\overline{y}_t^{\tau,v}|\mathcal{G}_t]=& \Big[\alpha_t\mathbb{E}[\overline{y}_t^{\tau,v}|\mathcal{G}_t]
+\beta_t\mathbb{E}[\overline{z}_t^{1,\tau,v}|\mathcal{G}_t]
+\mathbb{E}[\overline{z}_t^{1,\tau,v}|\mathcal{G}_t]
\cdot b(t,{\mathbb{E}[X_t|\mathcal{G}_t]},v_t)
+g(t,{\mathbb{E}[X_t|\mathcal{G}_t]},v_t)\Big]dt\\
&- \mathbb{E}[\overline{z}_t^{1,\tau,v}|\mathcal{G}_t]dW^1_t,\ t\in[0,\tau],\\
\mathbb{E}[\overline{y}_\tau^{\tau,v}|\mathcal{G}_\tau]=&\mathbb{E}\Big[\Phi(X_T)I_{\{\tau=T\}}
+{h}(\tau,X_\tau)I_{\{\tau<T\}}\Big|\mathcal{G}_\tau\Big].
\end{aligned}
\right.
\end{equation}
Using  the classical comparison theorem and the stable property of BSDE, we get the results as follows.
\begin{lemma}\label{th8}
Under the Assumption (H2), it holds, for each fixed $\tau\in\mathcal{T}_{0,T}$,
\begin{equation}\nonumber
\sup_{v\in\mathcal{V}}\overline{y}_0^{\tau,v}
=\widetilde{y}^\tau_0,
\end{equation}
 where $(\widetilde{y}^{\tau},\widetilde{z}^{1,\tau})$ is the unique $\mathbb{G}$-adapted solution of the BSDE
 \begin{equation}\nonumber
\left\{
\begin{aligned}
-d\widetilde{y}_t^{\tau}=&\Big[ \alpha_t\widetilde{y}_t^{\tau}+ {G}(t,{\mathbb{E}[X_t|\mathcal{G}_t]},
\widetilde{z}_t^{1,\tau})\Big]dt - \widetilde{z}_t^{1,\tau}dW^1_t,\ t\in[0,\tau],\\
\widetilde{y}_\tau^{\tau}=&
\mathbb{E}\Big[\Phi(X_T)I_{\{\tau=T\}}
+{h}(\tau,X_\tau)I_{\{\tau<T\}}\Big|\mathcal{G}_\tau\Big],
\end{aligned}
\right.
\end{equation}
where
$
G(t,x,z_1):=\beta_tz_1+\sup_{v\in U}\big\{z_1\cdot b(t,x,v)
+g(t,x,v)\big\}.
$
Moreover,
an optimal control has the following feedback form
$$v^*_t=\tilde{v}(t,\tilde{z}^{1,\tau}_t,{\mathbb{E}[X_t|\mathcal{G}_t]}),$$
where $\tilde{v}(t,z^1,x)=argmax_{v\in U}\big\{z_1\cdot b(t,x,v)
+g(t,x,v)\big\}.$
\end{lemma}
Combining  Lemma \ref{th8} and the classical result on  relationship between RBSDEs and optimal stopping problems (see, e.g., Theorem 3.3 in \cite{QS2014}), we obtain the characterization of the value of linear BRR problem \eqref{070806} as follows.
\begin{theorem}\label{th6}
Suppose that Assumption (H2) holds. Then it holds
\begin{equation}\nonumber
\sup_{v\in\mathcal{V}}{Y_0^v}(=\sup_{\tau\in\mathcal{T}_{0,T}}
\sup_{v\in\mathcal{V}}\overline{y}_0^{\tau,v}
=\sup_{\tau\in\mathcal{T}_{0,T}}\widetilde{y}_0^{\tau}
)=\mathcal{P}_0,
\end{equation}
where $(\mathcal{P},\mathcal{Q},\mathcal{K})$ is the solution of the following filtered RBSDE
\begin{equation}\nonumber
\left\{
\begin{aligned}
&\mathcal{P}_t= \mathbb{E}[\Phi(X_T)|\mathcal{G}_T]+\int_t^T\Big[\alpha_s\mathcal{P}_s+
{G}(s,{\mathbb{E}[X_s|\mathcal{G}_s]},
\mathcal{Q}_s)\Big]ds +\mathcal{K}_T-\mathcal{K}_t-\int_t^T\mathcal{Q}_sdW^1_s,\ t\in[0,T],\\
&\mathcal{P}_t-\mathbb{E}[h(t,X_t)|\mathcal{G}_t]\geq 0,\ \forall\ t\in[0,T],\ \text{a.s.},\ \int_0^T\mathcal{P}_{t}-\mathbb{E}[h(t,X_t)|\mathcal{G}_t]
d\mathcal{K}_t=0.
\end{aligned}
\right.
\end{equation}
Moreover, an optimal stopping is given as $\tau^*:=\inf\{t\geq 0: \mathcal{P}_t-\mathbb{E}[h(t,X_t)|\mathcal{G}_t]=0\}\wedge T.$
\end{theorem}

\subsection{Weak formulation of Convex BRR problems}
In this subsection, we  extend the above linear BRR problems to convex case. The goal is described by
\begin{equation}\label{071302}
{\bf (Convex\ BRR)}\ \sup_{v\in\mathcal{V}}Y^v_0,
\end{equation}
 where the payoff $Y_0^{v}$ is given by  the following convex conditional RBSDE
\begin{equation}\label{071301}
\left\{
\begin{aligned}
&Y_t^{v}=\Phi(X_T)+\int_t^Tf(s,{\mathbb{E}[X_s|\mathcal{G}_s]},Y_s^v,Z_s^{1,v},v_s)ds +K_T^{v}-K_t^{v}-\int_t^T(Z_s^{1,v},Z_s^{2,v})d{W^v_s},\ t\in[0,T],\\
&\mathbb{E}[Y_t^{v}- h(t,X_t)|\mathcal{G}_t]\geq 0,\ \forall\ t\in[0,T],\ \text{a.s.},\
\int_0^T\mathbb{E}[Y_t^{v}- h(t,X_t)|\mathcal{G}_t]dK_t^{v}=0,
\end{aligned}
\right.
\end{equation}
and $X$ is the solution of SDE \eqref{031801}, $\mathbb{P}^v$ and ${W}^v$ is given in \eqref{031802} and \eqref{032207}, respectively.
We assume that the terminal $\Phi$ and the barrier $h$ satisfy the Assumption (H2), and the driver $f$ satisfies \\
\begin{equation}\nonumber
{\bf (H3)}
\left\{
\begin{array}{l}
\text{(i)}\  f\ \text{is\ Lipschitz (with Lipschitz constant}\ \mu)\ \text{and\ {convex} in}\ (y,z),\ \text{uniformly\ in}\ (t,x,v).\\
\text{(ii)}\ f(s,0,0,0,0)\ \text{is}\ \mathbb{G}\text{-adapted};\ |f(t,x,0,0,v)|\leq C(1+|x|+|v|),\ \forall\ t\in[0,T].\\
 \end{array}
 \right.
 \end{equation}
 Notice that equation \eqref{071301} has a unique solution $(Y^v,Z^v,K^v)$ for each $v\in\mathcal{V}$, since it is equivalent to
\begin{equation}\nonumber
\left\{
\begin{aligned}
&Y_t^{v}=\Phi(X_T)+\int_t^T
\phi(s,{\mathbb{E}[X_s|\mathcal{G}_s]},Y_s^v,Z_s^{1,v},v_s)ds
+K_T^{v}-K_t^{v}-\int_t^T(Z_s^{1,v},Z_s^{2,v})d{W_s},\ t\in[0,T],\\
&\mathbb{E}[Y_t^{v}- h(t,X_t)|\mathcal{G}_t]\geq 0,\ \forall\ t\in[0,T],\ \text{a.s.},\
\int_0^T\mathbb{E}[Y_t^{v}- h(t,X_t)|\mathcal{G}_t]dK_t^{v}=0,
\end{aligned}
\right.
\end{equation}
where
{ $$\phi(t,x,y,z_1,v):=f(t,x,y,z_1,v)+
b(t,x,v)\cdot z_1.$$
It is easy to check that $\phi$ satisfies (H3).}
Using Fenchel-Moreau Theorem, we have
\begin{equation}\label{061502}
\phi(t,\mathbb{E}[X_t|\mathcal{G}_t],y,z^1,v_t)
=\esssup_{(\alpha,\beta)\in\mathcal{A}^\mathbb{G}}
\{\alpha_t y+ \beta_t z^1 -F(t,\mathbb{E}[X_t|\mathcal{G}_t],\alpha_t,\beta_t,v_t)\},
\end{equation}
where
 $\mathcal{A}^\mathbb{G}:=\big\{(\alpha,\beta):\ \mathbb{G}\text{-predictable},\ [-\mu,\mu]\times[-\mu,\mu] \text{-valued\ and}\ \mathbb{E}\int_0^T|F(t,\alpha_t,\beta_t,v_t)|^2dt
<\infty\big\}$ and
 \begin{equation}\nonumber
F(t,x,\alpha,\beta,v):=\sup_{(y,z^1)\in\mathbb{R}\times\mathbb{R}}
\{ \alpha y+ \beta z^1 -\phi(t,x,y,z^1,v)\}.
\end{equation}
We first establish the relationship of the solution of conditional RBSDEs between convex and linear drivers.
\begin{lemma}\label{le5.1}
For each $v\in\mathcal{V}$, it holds,
\begin{equation}\label{091703}
{Y^v_0}
={\sup_{
 (\alpha,\beta)\in\mathcal{A}^\mathbb{G}}Y^{v,\alpha,\beta}_0},
 \end{equation}
 where
 $(Y^{v,\alpha,\beta},Z^{v,\alpha,\beta},K^{v,\alpha,\beta})$
is the unique solution of the following linear conditional RBSDE
\begin{equation}\label{07170101}
\left\{
\begin{aligned}
&Y_t^{v,\alpha,\beta}
=\Phi(X_T)
+\int_t^T\Big[\alpha_sY_s^{v,\alpha,\beta}
+\beta_sZ_s^{1,v,\alpha,\beta}
-F(s,{\mathbb{E}[X_s|\mathcal{G}_s]},\alpha_s,\beta_s,v_s)\Big]ds\\
&\qquad\qquad +K_T^{v,\alpha,\beta}-K_t^{v,\alpha,\beta}
-\int_t^T(Z_s^{1,v,\alpha,\beta},Z_s^{2,v,\alpha,\beta})
d{W_s},\ t\in[0,T],\\
&\mathbb{E}[Y_t^{v,\alpha,\beta}- h(t,X_t)|\mathcal{G}_t]\geq 0,\ \forall\ t\in[0,T];\
\int_0^T\mathbb{E}[Y_t^{v,\alpha,\beta}- h(t,X_t)|\mathcal{G}_t]dK_t^{v,\alpha,\beta}=0.
\end{aligned}
\right.
\end{equation}
\end{lemma}
\begin{proof}
We denote $$c_s^1(\alpha,\beta):=\phi(s,{\mathbb{E}[X_s|\mathcal{G}_s]},Y_s^v,Z_s^{1,v},v_s)-\alpha_s Y_s^v-\beta_s Z_s^{1,v},\ c_s^2(\alpha,\beta):=-F(s,{\mathbb{E}[X_s|\mathcal{G}_s]},\alpha_s,\beta_s,v_s).$$
For each $(\alpha,\beta)\in\mathcal{A}^\mathbb{G}$, from \eqref{061502}  we have
$c_s^1(\alpha,\beta)\geq c_s^2(\alpha,\beta).$
Then using Corollary \ref{Co4.1} we obtain
 \begin{equation}\label{091701}
 Y^v_0
\geq Y^{v,\alpha,\beta}_0,\  \text{for\ each}\ (\alpha,\beta)\in\mathcal{A}^\mathbb{G}.
\end{equation}
On the other hand, from \eqref{061502} there exists $(\alpha^\varepsilon,\beta^\varepsilon)\in\mathcal{A}^\mathbb{G}$ (see, e.g. Lemma 3.1 in \cite{LL2019} for the construction of such $(\alpha^\varepsilon,\beta^\varepsilon)$)
such that
$$\phi(s,{\mathbb{E}[X_s|\mathcal{G}_s]},y,z^1,v_s)\leq \alpha_s^\varepsilon y+\beta_s^\varepsilon z^1-F(s,\mathbb{E}[X_s|\mathcal{G}_s],\alpha_s^\varepsilon,\beta_s^\varepsilon,v_s])
+\varepsilon.$$
Then using Theorem \ref{le4.1}, we get
\begin{equation}\nonumber
\begin{aligned}
\mathbb{E}\Big[\sup_{0\leq s\leq T}|Y_s^v-Y_s^{v,\alpha^\varepsilon,\beta^\varepsilon}|^2+\int_0^T
|Z_s^v-Z_s^{v,\alpha^\varepsilon,\beta^\varepsilon}|^2ds
+\sup_{0\leq s\leq T}|(K_T^v-K_s^v)-(K_T^{v,\alpha^\varepsilon,\beta^\varepsilon}-K_s^{v,\alpha^\varepsilon,\beta^\varepsilon})|^2\Big]
\leq C\varepsilon^2,
\end{aligned}
\end{equation}
from which we conclude that
$Y_0^v\leq Y_0^{v,\alpha^\varepsilon,\beta^\varepsilon}+C\varepsilon^{\frac12}$. Combining this and \eqref{091701}, we finally show \eqref{091703}.
\end{proof}
Using Lemma \ref{le5.1}, the convex BRR problem \eqref{071302} can be transformed to the supremum of a family of linear BRR problems
\begin{equation}\label{071501}
\begin{aligned}
 \sup_{v\in\mathcal{V}}{Y^v_0}
&=\sup_{v\in\mathcal{V}}{\sup_{
 (\alpha,\beta)\in\mathcal{A}^\mathbb{G}}Y^{v,\alpha,\beta}_0}=\sup_{
 (\alpha,\beta)\in\mathcal{A}^\mathbb{G}}
\Big(\sup_{v\in\mathcal{V}}Y^{v,\alpha,\beta}_0\Big).
\end{aligned}
\end{equation}
For each fixed $(\alpha,\beta)\in\mathcal{A}^\mathbb{G}$, we get from Corollary \ref{Co4.1} and Theorem \ref{le4.1} that
\begin{equation}\label{101603}
\sup_{v\in\mathcal{V}}Y^{v,\alpha,\beta}_0=
 \overline{Y}_0^{\alpha,\beta},
 \end{equation}
where $(\overline{Y}^{\alpha,\beta},\overline{Z}^{\alpha,\beta},\overline{K}^{\alpha,\beta})$
is the unique solution of the following linear conditional RBSDE
\begin{equation}\nonumber
\left\{
\begin{aligned}
&\overline{Y}_t^{\alpha,\beta}
=\Phi(X_T)
+\int_t^T\Big[\alpha_s\overline{Y}_s^{\alpha,\beta}
+\beta_s\overline{Z}_s^{1,\alpha,\beta}
-\overline{F}(s,{\mathbb{E}[X_s|\mathcal{G}_s]},\alpha_s,\beta_s)\Big]ds +\overline{K}_T^{\alpha,\beta}-\overline{K}_t^{\alpha,\beta}
\\
&\quad\quad\quad -\int_t^T(\overline{Z}_s^{1,\alpha,\beta},
\overline{Z}_s^{2,\alpha,\beta})dW_s,\ t\in[0,T],\\
&\mathbb{E}[\overline{Y}_t^{\alpha,\beta}- h(t,X_t)|\mathcal{G}_t]\geq 0,\ \forall\ t\in[0,T];\
\int_0^T\mathbb{E}[\overline{Y}_t^{\alpha,\beta}- h(t,X_t)|\mathcal{G}_t]d\overline{K}_t^{\alpha,\beta}=0,
\end{aligned}
\right.
\end{equation}
and
$\overline{F}(s,x,\alpha,\beta):=\inf_{v\in U}F(s,x,\alpha,\beta,v)$.
\begin{remark}
Compared with Theorem \ref{th6}, the proof of the conclusion \eqref{101603} seems more direct since the term involving in $z$ of the driver in \eqref{07170101} does not depend on the control $v$, different from the situation considered in \eqref{070805}. As a result, we can use the comparison theorem of conditional RBSDEs (Corollary \ref{Co4.1}) directly.
\end{remark}
Combining with \eqref{071501} and \eqref{101603}, similar to the proof of Lemma \ref{le5.1}, we obtain the following characterization of the value of convex BRR problem \eqref{071302} via the associated convex conditional RBSDE.
\begin{theorem}\label{th5.3}
The value of convex BRR problem \eqref{071302} has the representation
$$\sup_{v\in\mathcal{V}}Y^v_0=\overline{Y}_0,$$
where $(\overline{Y},\overline{Z},\overline{K})$ is the unique solution of the following convex conditional RBSDE
\begin{equation}\nonumber
\left\{
\begin{aligned}
&\overline{Y}_t
=\Phi(X_T)
+\int_t^T\overline{f}(s,{\mathbb{E}[X_s|\mathcal{G}_s]},\overline{Y}_s,\overline{Z}^1_s)ds +\overline{K}_T-\overline{K}_t
-\int_t^T(\overline{Z}_s^1,\overline{Z}_s^2)dW_s,\ t\in[0,T],\\
&\mathbb{E}[\overline{Y}_t- h(t,X_t)|\mathcal{G}_t]\geq 0,\ \forall\ t\in[0,T],\
\int_0^T\mathbb{E}[\overline{Y}_t- h(t,X_t)|\mathcal{G}_t]d\overline{K}_t=0,
\end{aligned}
\right.
\end{equation}
and the convex driver $\overline{f}$ is defined as follows
$$\overline{f}(t,{\mathbb{E}[X_t|\mathcal{G}_t]},y,z^1):=
\esssup_{(\alpha,\beta)\in\mathcal{A}^\mathbb{G}}
\{\alpha_t y+ \beta_tz^1 -\overline{F}(t,{\mathbb{E}[X_t|\mathcal{G}_t]},\alpha_t,\beta_t)\}.$$
\end{theorem}

\section{Backward recursive reflected control and zero-sum stochastic differential game problems with full information}


In contrast to the study of weak formulations of linear and convex control problems in Section 4, strong formulations of the general (requiring neither linear nor convex) BRR problems and a class of zero-sum stochastic differential games will be considered respectively in this section, but under  full information framework, i.e., $\mathbb{G}=\mathbb{F}$. Moreover, for both cases the state equations are driven by controlled stochastic functional differential equations, and both payoffs are described by the solution of the related RBSDEs.

For BRR problems, we show the value of the strong formulation is equal to that of weak ones. Such idea was firstly introduced by Bouchard, Elie, Moreau \cite{BEM2018} to address a type of linear control problems without any constraints on the recursive payoffs. Combining with nonlinear Snell envelope theory, we generalize the equivalent result between strong and weak formulations obtained in \cite{BEM2018} to a type of nonlinear control problems with constrained payoffs.
Then we characterize the value of strong formulation of BRR problems via the associated RBSDEs.
Moreover, we extend the study of the general BRR problems to a type of zero-sum stochastic differential games and obtain the closed form of the saddle point with the help of the solution of the corresponding RBSDE under the well-known Isaacs condition.

In this section, the underlying probability space
$(\Omega,{\mathcal{F}},\mathbb{P})$ is chosen to be a Wiener space, namely, $\Omega=C_0([0,T];\mathbb{R}^d)$ is the set of all continuous functions from $[0,T]$ to $\mathbb{R}^d$ with value $0$ at initial time, $\mathcal{F}$ is the complete Borel $\sigma$-field on $\Omega$, $\mathbb{P}$ is the Wiener measure such that the canonical processes $W_s(\omega)=\omega(s)$, $s\in[0,T]$, $\omega\in\Omega$, is a $d$-dimensional standard Brownian motion.
\subsection{Strong formulation of BRR problems with full information}
We formulate the strong version of BRR problems with full information. In this situation,
the set $\mathcal{V}$ of admissible controls in Section 4  turns out to be
\begin{equation}\label{110601}
\mathcal{V}_{\mathbb{F}}:=\Big\{v\Big|v\ \text{is}\ V\text{-valued}\ \mathbb{F}\text{-adapted\ process\ such\ that}\ \mathbb{E}\Big[\int_0^T|v_t|^2dt\Big]<\infty\Big\}.
\end{equation}
We denote by $\mathbb{X}$ the space of continuous functions from $[0,T]$ to $\mathbb{R}^d$ endowed with the uniform norm $\|X\|_t=\sup_{0\leq s\leq t}|X_t|$, $t\in[0,T]$.
Let the measurable functions
$b:[0,T]\times\mathbb{X}\times V\rightarrow\mathbb{R}^d$ and $\sigma:[0,T]\times\mathbb{X}\rightarrow\mathbb{R}^{d\times d}
$
 satisfy
 \begin{equation}\nonumber
{\bf (A1)}
\left\{
\begin{array}{l}
\text{(i)}\ \text{For\ each}\ v\in{V}\ \text{and\ continuous\
progressively\ measurable\ process}\ x, (b(t,x,v))_{0\leq t\leq T}\ \text{and}\\
\quad\  (\sigma(t,x))_{0\leq t\leq T}\ \text{are}\ \text{progressively\ measurable.}\\
\text{(ii)}\ \text{There\ exists\ a\ constant}\ C>0\ \text{such\ that},\ \text{for\ every}\ t\in[0,T],\ x,x'\in\mathbb{X},\ v\in V,\\
\qquad\qquad |b(t,x,v)-b(t,x',v)|+|\sigma(t,x)-\sigma(t,x')|\leq C\|x-x'\|_{t}.\\
\text{(iii)}\ b\ \text{is\ continuous\ in}\ v\ \text{and\ uniformly\ bounded}.\\
\text{(iv)}\  \sigma\ \text{is\ uniformly\ bounded\ and invertible}, \text{and\ its\ inverse}\, \sigma^{-1}\ \text{is\ also\ uniformly\ bounded}.
 \end{array}
 \right.
 \end{equation}
The controlled state is described by the following stochastic functional differential equation (SFDE)
\begin{equation}\label{102601}
\begin{aligned}
X_t^{v}=x_0+\int_0^tb(s,X^{v},v_s)ds+\int_0^t\sigma(s,X^{v})dW_s,\ t\in[0,T].
\end{aligned}
\end{equation}
Under the condition (A1), it is well known (see, e.g., Theorem 2.1 in \cite{RY2013} (Page 348 therein)) that SFDE \eqref{102601} exists a unique solution $X^v\in \mathcal{S}^2$ for each $v\in\mathcal{V}_\mathbb{F}$.
\begin{remark}
Noting that the coefficients $b$ and $\sigma$ at each time $t$ depend on the entire path of the state $X^v$ from $0$ to $t$ rather than only the current time  $t$, such SFDE is non-Markovian.

On the other hand, the boundedness assumption of $b$, $\sigma$ and $\sigma^{-1}$ can be relaxed to the linear growth condition, such as $|b(t,x,v)|\leq C(1+\|x\|_t)$, $(t,x,v)\in[0,T]\times\mathbb{X}\times V$. We impose those stronger assumption as in (A1) in order to  avoid more technique details and focus on the novelty of our approach.
\end{remark}

We introduce the associated constrained recursive payoff. Let
$f:[0,T]\times\mathbb{X}\times\mathbb{R}\times\mathbb{R}^d\times V\rightarrow\mathbb{R},\ h:[0,T]\times\mathbb{X}\rightarrow\mathbb{R},\ \Phi:\mathbb{X}\rightarrow\mathbb{R},\ \text{satisfy}$
 \begin{equation}\nonumber
{\bf (A2)}
\left\{
\begin{array}{l}
\text{(i)}\ \text{For\ each}\ (y,z,v)\in\mathbb{R}\times\mathbb{R}^d\times V\
\text{and\ continuous\
progressively\ measurable\ process}\ x,\\
\quad\ f(\cdot,x,y,z,v),\ h(\cdot,x)\ \text{and}\ \Phi(x)\ \text{are}\ \text{progressively\ measurable};\\
\text{(ii)}\ f\  \text{is\ {continuous}\ in}\ (t,v)\  \text{and\ there\ exists\ a\ constant}\
C>0\ \text{such\ that},\ \text{for\ all}\ t\in[0,T],\ v\in  V, \\
\quad\ x,x'\in\mathbb{X},\ (y,z),(y',z')\in\mathbb{R}\times\mathbb{R}^d,\ \\
\qquad\qquad\quad |f(t,x,y,z,v)-f(t,x',y',z',v)|\leq C(\|x-x'\|_t+|y-y'|+|z-z'|);\\ 
\text{(iii)}\ h\ \text{is\ continuous in}\ {(t,x)}\ \text{and\ there exists a\ constant}\ C>0\ \text{such\ that}\
{|h(t,x)|\leq C(1+\|x\|_t)}; \\
\text{(iv)}\ 
\text{There\ exists\ a\ constant}\
C>0\ \text{such\ that,}\ {|\Phi(x)|\leq C(1+\|x\|_T)};\ h(T,x)\leq \Phi(x),\ x\in\mathbb{X}.
 \end{array}
 \right.
 \end{equation}
The constrained  payoff $Y_0^v$ with the admissible control $v$ is described
 by  the following controlled RBSDE
\begin{equation}\nonumber
\left\{
\begin{aligned}
&Y_t^{v}=\Phi(X^{v})+\int_t^Tf(s,X^{v},Y_s^{v},Z_s^{v},v_s)ds +K_T^{v}-K_t^{v}-\int_t^TZ_s^{v}dW_s,\ t\in[0,T],\\
&Y_t^{v}\geq h(t,X^{v}),\ \forall\ t\in[0,T],\ \text{a.s.},\ \int_0^T[Y_t^{v}- h(t,X^{v})]dK_t^{v}=0,
\end{aligned}
\right.
\end{equation}
where $\Phi$ and $f$ stands for the terminal and instantaneous payoff, respectively, $h$ is the constraint condition of the payoff.
The functionals $\Phi,f$ and $h$ are allowed to rely on the entire history state rather than only the current value.
For each $v\in\mathcal{V}_\mathbb{F}$, it is clear that there exists a unique solution $(Y^v,Z^v,K^v)\in \mathcal{S}^2\times\mathcal{H}^2\times\mathcal{A}^2$ under the condition (A2).
The aim of the controller is to maximize this payoff $Y_0^v$ over all admissible controls, i.e.,
\begin{equation}\label{091704}
{ \textbf{(Strong\ BRR-F)}}\sup_{v\in\mathcal{V}_\mathbb{F}}Y_0^{v}.
\end{equation}

In order to address the strong BRR-F problem \eqref{091704}, we introduce the weak formulation of this problem and then show their values coincide. We denote by $X$ the unique solution of the following SFDE
\begin{equation}\label{091718}
\begin{aligned}
X_t=x_0+\int_0^t \sigma(s,X)dW_s,\ t\in[0,T].
\end{aligned}
\end{equation}
It is clear that $
\mathbb{E}[\|X\|_T^p]\leq C(1+|x_0|^p),\ \text{for\ all}\ p\geq 2.
$
For each given admissible control  $v\in\mathcal{V}_\mathbb{F}$, we define a probability measure $\mathbb{P}^{v}$ on $(\Omega,\mathcal{F})$, which is equivalent to $\mathbb{P}$ and whose density function is given by
\begin{equation}\nonumber
\frac{d\mathbb{P}^{v}}{d\mathbb{P}}\Big|_{\mathcal{F}_T}
=\exp\{\int_0^T\sigma^{-1}(t,X)b(t,X,v_t)dW_t
-\frac12\int_0^T|\sigma^{-1}(t,X)b(t,X,v_t)|^2dt\}.
\end{equation}
Thanks to Girsanov Theorem, the process
\begin{equation}\nonumber
dW^{v}_t:=-\sigma^{-1}(t,X)b(t,X,v_t)dt+dW_t,\ t\in[0,T],
 \end{equation}
is a Brownian motion under the probability measure $\mathbb{P}^{v}$. Moreover, $X$ is the weak solution of the following SFDE
\begin{equation}\nonumber
\begin{aligned}
X_t=x_0+\int_0^t b(s,X,v_s)ds+\int_0^t\sigma(s,X)dW_s^{v},\ t\in[0,T].
\end{aligned}
\end{equation}
The aim of this weak formulation of BRR-F problem is given by
\begin{equation}\label{091716}
{ \textbf{(Weak\ BRR-F)}}\sup_{v\in\mathcal{V}_\mathbb{F}}\mathcal{Y}_0^{v},
\end{equation}
where $(\mathcal{Y}^{v},\mathcal{Z}^{v},\mathcal{K}^{v})$ is the solution of the following controlled RBSDE
\begin{equation}\nonumber
\left\{
\begin{aligned}
&\mathcal{Y}_t^{v}=\Phi(X)+\int_t^Tf(s,X,\mathcal{Y}_s^{v},\mathcal{Z}_s^{v},v_s)ds +\mathcal{K}_T^{v}-\mathcal{K}_t^{v}-\int_t^T\mathcal{Z}_s^{v}dW_s^{v},\ t\in[0,T],\\
&\mathcal{Y}_t^{v}\geq h(t,X),\ \forall\ t\in[0,T],\ \int_0^T[\mathcal{Y}_t^{v}- h(t,X)]d\mathcal{K}_t^{v}=0.
\end{aligned}
\right.
\end{equation}
Then we have the following relationship between the strong BRR-F problem \eqref{091704} and weak ones \eqref{091716}.
\begin{theorem}\label{th6.01}
Under the Assumptions (A1)-(A2), it holds
\begin{equation}\nonumber
\sup_{v\in\mathcal{V}_\mathbb{F}}\mathcal{Y}_0^{v}
=\sup_{v\in\mathcal{V}_\mathbb{F}}{Y}_0^{v}.
\end{equation}
\end{theorem}
\begin{proof}
We denote by
$\mathcal{T}_{t,T}^\mathbb{F}$ the set of $\mathbb{F}$-stopping times with values in $[t,T]$. It follows from  the nonlinear Snell envelope theory (see, e.g., Theorem 3.3 in \cite{QS2014}),  for each $v\in\mathcal{V}_\mathbb{F}$, $t\in[0,T]$,
\begin{equation}\label{091705}
Y_t^{v}=\esssup_{\tau\in\mathcal{T}^\mathbb{F}_{t,T}}
y_t^{\tau,v},
\end{equation}
where, for each $\tau\in\mathcal{T}^\mathbb{F}_{t,T}$, $(y^{\tau,v},z^{\tau,v})$ is the unique solution of the following BSDE
\begin{equation}\nonumber
\left.
\begin{aligned}
y_s^{\tau,v}=\Big[\Phi(X^{v})I_{\{\tau=T\}}+{h}(\tau,X^{v})I_{\{\tau<T\}}\Big]
+\int_s^\tau f(r,X^{v},y_r^{\tau,v},z_r^{\tau,v},v_r)dr - \int_s^\tau z_r^{\tau,v}dW_r,\ s\in[t,\tau].
\end{aligned}
\right.
\end{equation}
From \eqref{091705}, we have
\begin{equation}\label{091714}
Y_0^{v}=\sup_{\tau\in\mathcal{T}^\mathbb{F}_{0,T}}
y_0^{\tau,v},\ \text{and\ similarly}\ \mathcal{Y}_0^{v}=\sup_{\tau\in\mathcal{T}^\mathbb{F}_{0,T}}
\overline{y}_0^{\tau,v},
\end{equation}
where, for each $\tau\in\mathcal{T}^\mathbb{F}_{0,T}$, $(\overline{y}^{\tau,v},\overline{z}^{\tau,v})$ is the unique solution of the following BSDE
 \begin{equation}\nonumber
\overline{y}_t^{\tau,v}=\Big[\Phi(X)I_{\{\tau=T\}}+{h}(\tau,X)
I_{\{\tau<T\}}\Big]+\int_t^\tau f(s,X,\overline{y}_s^{\tau,v},\overline{z}_s^{\tau,v},v_s)ds - \int_t^\tau\overline{z}_s^{\tau,v}dW_s^{v},\ t\in[0,\tau].
\end{equation}
Step 1. We show that for each $(\tau,v)\in\mathcal{T}^\mathbb{F}_{0,T}\times\mathcal{V}_1$, there exist $(\tau_1,v_1), (\tau_2,v_2)\in\mathcal{T}^\mathbb{F}_{0,T}
\times\mathcal{V}_\mathbb{F}$ such that
\begin{equation}\label{091709}
y_0^{\tau,v}=\overline{y}_0^{\tau_1,v_1},\ \overline{y}_0^{\tau,v}=y_0^{\tau_2,v_2},
\end{equation}
where $\mathcal{V}_1\subseteq\mathcal{V}_\mathbb{F}$ is the set of simple processes $v$, i.e.,
\begin{equation}\label{091706}
 v(t,\omega)=\sum_{i=0}^{N-1}\zeta_i(\omega)\cdot I_{\{t_{i}< t\leq t_{i+1}\}},
\end{equation}
where $\pi=\{0=t_0<t_1<\cdots<t_N=T\}$ is a partition of $[0,T]$, $\zeta_i$ is $\mathcal{F}_{t_i}$-measurable bounded $V$-valued random variable, $i=0,1,2,\cdots,N-1.$

For each $v\in\mathcal{V}_1$ with the form \eqref{091706}, we identify $\zeta_i$ as a Borel measurable function $\omega\rightarrow \zeta_i(\omega)=\zeta_i(\omega_{\cdot\wedge t_i})$, and we define
$${v}_1(t,\omega)=\sum_{i=0}^{N-1}\zeta_i(\omega^{\zeta})\cdot I_{\{t_{i}< t\leq t_{i+1}\}},$$
where $w^{\zeta}$ is defined recursively as follows, for $i=0,1,2,\cdots,N-1,$ $t\in (t_i,t_{i+1}]$,
\begin{equation}\label{091707}
\omega^{\zeta}_0=0,\ \omega^{\zeta}_t=\omega_t-\sum_{k=0}^{i-1}\int_{t_k}^{t_{k+1}}
\sigma^{-1}(s,X)b(s,X,\zeta_k(\omega_{t_k}^\zeta))ds-
\int_{t_i}^{t}
\sigma^{-1}(s,X)b(s,X,\zeta_i(\omega_{t_i}^\zeta))ds.
\end{equation}
It is easy to check that ${v}_1\in\mathcal{V}_\mathbb{F}$.
Comparing the following two SFDEs
\begin{equation}\nonumber
\left\{
\begin{aligned}
X_t^{v}=&X_{t_i}^{v}+\int_{t_i}^tb\big(s,X^{v},\zeta_i(W_{s\wedge t_i})\big)ds+\int_{t_i}^t\sigma(s,X^{v})dW_s,\ t\in[t_i,t_{i+1}],\\
X_t=&X_{t_i}+\int_{t_i}^tb\big(s,X,\zeta_i(W_{s\wedge t_i}^{{v_1}})\big)ds+\int_{t_i}^t\sigma(s,X)dW_s^{{v_1}},\ t\in[t_i,t_{i+1}],
\end{aligned}
\right.
\end{equation}
we obtain  the law of $(X^{v},v,W)$ under $\mathbb{P}$ and that of $(X,{v}_1,W^{{v_1}})$ under $\mathbb{P}^{v_1}$ coincide from the uniqueness of the weak solution of SFDE (see, e.g., Theorem 4.2 of Chapter 4 in \cite{IW2014}). For each $\tau\in\mathcal{T}^\mathbb{F}_{0,T}$, we define
\begin{equation}\label{091710}
\tau_1(\omega):=\tau(\omega^{\zeta}),
\end{equation}
where $\omega^{\zeta}$ is given in \eqref{091707}. {Since $\omega\rightarrow\tau(\omega)=\tau(W_\cdot)$ can be identified as a Borel measurable function, then
$\tau_1(\omega)=\tau(W^{v_1}_\cdot)$
 is a stopping time, i.e., $\tau_1\in\mathcal{T}^\mathbb{F}_{0,T}$.}
Using the discrete-time approximation (see, e.g., Lemma A.4 in \cite{BEM2018}) for the following BSDEs with $\widetilde{h}(t,x):=\Phi(x)I_{\{t=T\}}+h(t,x)I_{\{t<T\}}$,
\begin{equation}\nonumber
\begin{aligned}
y_0^{\tau,v}=&\widetilde{h}(\tau,X^{v})+\int_0^T
I_{\{t\leq\tau\}}\cdot f(t,X^{v},y_t^{\tau,{v}},z_t^{\tau,{v}},v(t))dt
-\int_0^TI_{\{t\leq\tau\}}\cdot z_t^{\tau,v}dW_t,\\
\overline{y}_0^{\tau_1,v_1}=&
\widetilde{h}(\tau_1,X)
+\int_0^TI_{\{t\leq\tau_1\}}\cdot f(t,X,\overline{y}_t^{\tau_1,{v_1}},\overline{z}_t^{\tau_1,{v_1}},v_1(t))dt
-\int_0^TI_{\{t\leq\tau_1\}}\cdot \overline{z}_t^{\tau_1,v_1}dW^{v_1}_t,
\end{aligned}
\end{equation}
we get
\begin{equation}\label{091708}
y_0^{\tau,v}=\lim_{n\rightarrow\infty}y_0^{n,\tau,v},\  \overline{y}_0^{\tau_1,{v_1}}=\lim_{n\rightarrow\infty}\overline{y}_0^{n,\tau_1,{v_1}},
\end{equation}
where $(y^{n,\tau,v},z^{n,\tau,v})$ and $(\overline{y}^{n,\tau_1,v_1},\overline{z}^{n,\tau_1,v_1})$ is defined recursively, respectively, as follows, for $i=n-1,\cdots,0$ (with $t_i^n:=i\frac{T}{n}$)
\begin{equation}\nonumber
\begin{aligned}
&y_{t_i^n}^{n,\tau,v}
=\mathbb{E}[y_{t_{i+1}^n}^{n,\tau,v}+\int_{t_i^n}^{t_{i+1}^n}
I_{\{t\leq\tau(W)\}}\cdot f(t,X^{v},{y}_{t_i^n}^{n,\tau,v},{z}_{t_i^n}^{n,\tau,v},\zeta_i(W_{t\wedge t_i^n}))dt
\big|\mathcal{F}_{t_i^n}],\\
& z_{t_i^n}^{n,\tau,v}=(t_{i+1}^n-t_{i}^n)^{-1}\cdot
I_{\{t_i^n\leq\tau(W)\}}\cdot
\mathbb{E}[y_{t_{i+1}^n}^{n,\tau,v}(W_{t_{i+1}^n}-W_{t_{i}^n})
\big|\mathcal{F}_{t_i^n}],\\
&\overline{y}_{t_i^n}^{n,\tau_1,v_1}
=\mathbb{E}_{\mathbb{P}^{v_1}}[\overline{y}_{t_{i+1}^n}^{n,\tau_1,v_1}
+\int_{t_i^n}^{t_{i+1}^n}
I_{\{t\leq\tau(W^{v_1})\}}\cdot f(t,X,\overline{y}_{t_i^n}^{n,\tau_1,v_1},
\overline{z}_{t_i^n}^{n,\tau_1,v_1},\zeta_i(W_{t\wedge t_i^n}^{v_1}))dt
\big|\mathcal{F}_{t_i^n}],\\
& \overline{z}_{t_i^n}^{n,\tau_1,v_1}=(t_{i+1}^n-t_{i}^n)^{-1}\cdot
I_{\{t_i^n\leq\tau(W^{v_1})\}}\cdot
\mathbb{E}_{\mathbb{P}^{v_1}}[\overline{y}_{t_{i+1}^n}^{n,\tau_1,v_1}
(W^{v_1}_{t_{i+1}^n}-W^{v_1}_{t_{i}^n})
\big|\mathcal{F}_{t_i^n}].
\end{aligned}
\end{equation}
Noting that the law of $(X^{v},W)$ under the probability measure $\mathbb{P}$ is the same to that of $(X,W^{{v}_1})$ under the measure $\mathbb{P}^{v_1}$, we obtain
$y_0^{n,\tau,v}=\overline{y}_0^{n,\tau,{v_1}}.$ Then it follows from \eqref{091708} that the first equality in \eqref{091709} holds, i.e.,  $y_0^{\tau,v}=\overline{y}_0^{\tau_1,{v_1}}.$

Similarly, we can show that the second equality in \eqref{091709} holds. In this case, for each $v\in\mathcal{V}_1$ with the form \eqref{091706}, we define
$$ {v}_2(s,\omega)=\zeta_i(\omega^{\zeta}_{s\wedge t_i}),\ s\in (t_i,t_{i+1}],$$
where $w^{\zeta}$ is defined recursively as follows: for $i=0,1,2,\cdots,n-1,$ $\omega^{\zeta}_0=0,$
\begin{equation}\nonumber
\begin{aligned}
\omega^{\zeta}_s=\omega_s+\sum_{k=0}^{i-1}\int_{t_k}^{t_{k+1}}
\sigma^{-1}(t,X^{v_2})b(t,X^{v_2},
\zeta_k(\omega_{t_k}^{\zeta}))dt
+
\int_{t_i}^{s}
\sigma^{-1}(t,X^{v_2})b(t,X^{v_2},
\zeta_i(\omega_{t_i}^{\zeta}))dt,\ s\in (t_i,t_{i+1}].
\end{aligned}
\end{equation}
For each $\tau\in\mathcal{T}_{0,T}^\mathbb{F}$, we define $\tau_2\in\mathcal{T}_{0,T}^\mathbb{F}$
similar to the definition of $\tau_1$ given in \eqref{091710}. {Using the same arguments as above}, it holds $\overline{y}_0^{\tau,v}=y_0^{\tau_2,{v}_2}$.
\\
Step 2. We show that
\begin{equation}\label{091711}
\sup_{v\in\mathcal{V}_\mathbb{F}}\sup_{\tau\in\mathcal{T}_{0,T}^\mathbb{F}}y_0^{\tau,v}
=\sup_{v\in\mathcal{V}_\mathbb{F}}
\sup_{\tau\in\mathcal{T}_{0,T}^\mathbb{F}}
\overline{y}_0^{\tau,v}.
\end{equation}
For each $v\in\mathcal{V}_\mathbb{F}$, there exists a sequence  $v^n\in\mathcal{V}_1$ such that $\mathbb{E}\int_0^T|v_s-v_s^n|^2ds\rightarrow 0,\ \text{as}\ n\rightarrow\infty.$ Then from the classical arguments, we get
$\mathbb{E}[\|X^{v}-X^{v^n}\|_T^2]\rightarrow 0.$ 
Moreover, from the stability property of BSDE, we get {$y_0^{\tau,v}=\lim_{n\rightarrow\infty}y_0^{\tau,v^n}$}. 
Then it holds
\begin{equation}\label{091713}
\sup_{v\in\mathcal{V}_\mathbb{F}}\sup_{\tau\in\mathcal{T}_{0,T}^\mathbb{F}}y_0^{\tau,v}
=
\sup_{v\in\mathcal{V}_1}\sup_{\tau\in\mathcal{T}_{0,T}^\mathbb{F}}y_0^{\tau,v}.
\end{equation}
Similarly, we have
\begin{equation}\label{091712}
\sup_{v\in\mathcal{V}_\mathbb{F}}\sup_{\tau\in\mathcal{T}_{0,T}^\mathbb{F}}\overline{y}_0^{\tau,v}
=
\sup_{v\in\mathcal{V}_1}\sup_{\tau\in\mathcal{T}_{0,T}^\mathbb{F}}\overline{y}_0^{\tau,v}.
\end{equation}
Using the result of Step 1 (i.e., \eqref{091709}), \eqref{091713} and \eqref{091712}, we obtain \eqref{091711}.

 Finally, combining \eqref{091714} and \eqref{091711}, {we get the desired result}.
\end{proof}

From Theorem \ref{th6.01}, we address the strong BRR-F problem \eqref{091704} via the weak BRR-F  problem \eqref{091716}.
For this,
we introduce the following Hamiltonian functional
\begin{equation}\nonumber
{ F(t,x,y,z,v):=f(t,x,y,z,v)
+z\sigma^{-1}(t,x)b(t,x,v)},\ (t,x,y,z,v)\in [0,T]\times\mathbb{X}\times\mathbb{R}\times\mathbb{R}^d\times V.
\end{equation}
Under the Assumptions (A1)-(A2), $F$ is Lipschitz in $(y,z)$, uniformly
with respect to $(t,x,v)$ and there exists a constant $C>0$ (independent of $v$) such that
\begin{equation}\nonumber
|F(t,x,y,z,v)|\leq C(1+\|x\|_t+|y|+|z|).
\end{equation}
We denote
$$G(t,x,y,z)=\sup_{v\in V}F(t,x,y,z,v),\ (t,x,y,z)\in[0,T]\times\mathbb{X}\times\mathbb{R}\times\mathbb{R}^d.$$
Since $F$ is continuous on the compact space $V$, there exists a measurable mappings
 $\bar{v}:[0,T]\times\mathbb{X}\times\mathbb{R}\times
\mathbb{R}^d\rightarrow V$ such that
\begin{equation}\label{091717}
G(t,x,y,z)
=F(t,x,y,z,\bar{v}(t,x,y,z)).
\end{equation}
Then using  comparison theorem of RBSDEs (see, for example, Theorem 4.1 in \cite{KKPPQ}) and Theorem \ref{th6.01}, we get the results as follows.
\begin{theorem}\label{th6.3}
Suppose that the Assumptions (A1)-(A2) hold. Then the value of the strong BRR-F problem \eqref{091704} can be characterized as follows
\begin{equation}\nonumber
\sup_{v\in\mathcal{V}_{\mathbb{F}}}
Y_0^v=
p_0,
\end{equation}
where $(p,q,k)$ is the solution of the following RBSDE
\begin{equation}\nonumber
\left\{
\begin{aligned}
&p_t= \Phi(X)+\int_t^TG(s,X,p_s,q_s)ds +k_T-k_t-\int_t^Tq_sdW_s,\ t\in[0,T],\\
&p_t\geq h(t,X),\ \forall\ t\in[0,T];\ \int_0^T[p_{t}-h(t,X)]dk_t=0.
\end{aligned}
\right.
\end{equation}
Moreover, an optimal control $v^*\in
\mathcal{V}_{\mathbb{F}}$ has the following feedback form
$$ v^*_t=\bar{v}(t,X,p_t,q_t),
$$
where the function $\bar{v}$ is given in \eqref{091717}.
\end{theorem}
\begin{remark}
When strong BRR-F problem \eqref{091704} is of \emph{Markovian type}, namely, all the involving coefficients $b,\sigma,\Phi,f$ and $h$ (at time $t$) rely on $X_t^v$ rather than $(X_s^v)_{0\leq s\leq t}$, such optimal control problem has been studied by Wu and Yu \cite{WY2008}   by using dynamic programming principle  approach. Compared with their work, the advantage of our approach is that it can be applied to address such strong BRR-F problems within non-Markovian framework. On the other hand, our approach can be applied directly to solve  zero-sum stochastic differential games as shown in the next subsection.
\end{remark}

\subsection{Zero-sum stochastic differential games with full information}
In this subsection, we generalize the strong BRR-F problem \eqref{091704} to  zero-sum stochastic differential game case. For this, let $U$  be a nonempty {compact} subset of $\mathbb{R}^{m}$. The admissible control space for Player 1 is denoted by $\mathcal{U}_\mathbb{F}$, which is defined similarly to the admissible control space $\mathcal{V}_\mathbb{F}$ (see \eqref{110601}) for Player 2
with  $V$ replacing by $U$.
We formulate the model of the game problem.
The controlled state is driven by the following SFDE
$$
X_t^{u,v}=x_0+\int_0^tb(s,X^{u,v},u_s,v_s)ds+\int_0^t\sigma(s,X^{u,v})dW_s,\ t\in[0,T].
$$
The payoff $J(u,v)$ is defined by
\begin{equation}\label{090801}
J(u,v)=Y_0^{u,v},
\end{equation}
where $(Y^{u,v},Z^{u,v},K^{u,v})$ is the solution of  the following controlled RBSDE
\begin{equation}\nonumber
\left\{
\begin{aligned}
&Y_t^{u,v}=\Phi(X^{u,v})+\int_t^Tf(s,X^{u,v},Y_s^{u,v},Z_s^{u,v},u_s,v_s)ds +K_T^{u,v}-K_t^{u,v}-\int_t^TZ_s^{u,v}dW_s,\ t\in[0,T],\\
&Y_t^{u,v}\geq h(t,X^{u,v}),\ \forall\ t\in[0,T];\ \int_0^T[Y_t^{u,v}- h(t,X^{u,v})]dK_t^{u,v}=0.
\end{aligned}
\right.
\end{equation}
Herein, $J$ represents the cost for Player 1 and the gain for Player 2. Thus, Player 1 aims to minimize  $J(u,v)$ by using the control $u$, while Player 2 wants to maximize $J(u,v)$ via the control $v$. For such zero-sum games,
we want  to find a saddle point $(u^*,v^*)\in\mathcal{U}_\mathbb{F}\times\mathcal{V}_\mathbb{F}$, i.e.,   for   all admissible control pair $(u,v)\in\mathcal{U}_\mathbb{F}\times\mathcal{V}_\mathbb{F}$, it holds
\begin{equation}\label{080203}
J(u^*,v)\leq J(u^*,v^*)\leq J(u,v^*).
\end{equation}
The coefficients $b,\sigma$ and $f,\Phi, h$
satisfy the same conditions  of those (A1) and (A2) in Subsection 5.1 with
the variable $v$ replacing by a pair of variables $(u,v)$. It is clear that the above
SFDE and RBSDE exist a unique solution $(X^{u,v},Y^{u,v},Z^{u,v},K^{u,v})$ for each admissible control pair
$(u,v)\in\mathcal{U}_\mathbb{F}\times\mathcal{V}_\mathbb{F}$.

In order to find the saddle point of \eqref{090801}, we introduce an auxiliary weak formulation of this game problem and then show that its saddle  point exists, which is also a saddle point for original problem \eqref{090801}.
The state equation of the auxiliary game problem is still described by  SFDE \eqref{091718}.
For each given admissible control pair $(u,v)\in\mathcal{U}_\mathbb{F}\times\mathcal{V}_\mathbb{F}$, we define an equivalent probability measure $\mathbb{P}^{u,v}$ on $(\Omega,\mathcal{F})$:
\begin{equation}\nonumber
\frac{d\mathbb{P}^{u,v}}{d\mathbb{P}}\Big|_{\mathcal{F}_T}
=\exp\{\int_0^T\sigma^{-1}(t,X)b(t,X,u_t,v_t)dW_t
-\frac12\int_0^T|\sigma^{-1}(t,X)b(t,X,u_t,v_t)|^2dt\}.
\end{equation}
 Then the process
$
W^{u,v}_t:=-\int_0^t\sigma^{-1}(s,X)b(s,X,u_s,v_s)ds+W_t,\ t\in[0,T],
$
is a Brownian motion under the probability measure $\mathbb{P}^{u,v}$ . 
The payoff $\mathcal{J}(u,v)$ of the auxiliary game problem is given by
\begin{equation}\label{09080101}
\mathcal{J}(u,v)=\mathcal{Y}_0^{u,v},
\end{equation}
where $(\mathcal{Y}^{u,v},\mathcal{Z}^{u,v},\mathcal{K}^{u,v})$ is the solution of the following controlled RBSDE
\begin{equation}\nonumber
\left\{
\begin{aligned}
&\mathcal{Y}_t^{u,v}=\Phi(X)+\int_t^Tf(s,X,\mathcal{Y}_s^{u,v},\mathcal{Z}_s^{u,v},u_s,v_s)ds +\mathcal{K}_T^{u,v}-\mathcal{K}_t^{u,v}-\int_t^T\mathcal{Z}_s^{u,v}dW_s^{u,v},\ t\in[0,T],\\
&\mathcal{Y}_t^{u,v}\geq h(t,X),\ \forall\ t\in[0,T];\ \int_0^T[\mathcal{Y}_t^{u,v}- h(t,X)]d\mathcal{K}_t^{u,v}=0.
\end{aligned}
\right.
\end{equation}

We have the following equivalent relation for these two game problems.
\begin{theorem}\label{th6.1}
The upper values (resp., the lower values) of game problems \eqref{090801} and \eqref{09080101} coincide, i.e.,
\begin{equation}\nonumber
\inf_{u\in\mathcal{U}_\mathbb{F}}\sup_{v\in\mathcal{V}_\mathbb{F}}J(u,v)
=\inf_{u\in\mathcal{U}_\mathbb{F}}\sup_{v\in\mathcal{V}_\mathbb{F}}\mathcal{J}(u,v),\
\sup_{v\in\mathcal{V}_\mathbb{F}}\inf_{u\in\mathcal{U}_\mathbb{F}}J(u,v)
=\sup_{v\in\mathcal{V}_\mathbb{F}}\inf_{u\in\mathcal{U}_\mathbb{F}}\mathcal{J}(u,v).
\end{equation}
\end{theorem}
We omit its proof since it  is similar to that of Theorem \ref{th6.01}.

From Theorem \ref{th6.1}, we address the original game problem \eqref{090801} via the equivalent auxiliary game problem \eqref{09080101}.
For this,
we introduce the following Hamiltonian functional
\begin{equation}\nonumber
{ F(t,x,y,z,u,v)=f(t,x,y,z,u,v)
+z\sigma^{-1}(t,x)b(t,x,u,v)},\ (t,x,y,z,u,v)\in[0,T]\times\mathbb{X}\times\mathbb{R}\times
\mathbb{R}^d\times U\times V.
\end{equation}
Similar to most researches on stochastic differential games, we assume that the following Isaacs condition holds:
$$G(t,x,y,z):=\inf_{u\in U}\sup_{v\in V}F(t,x,y,z,u,v)=\sup_{v\in V}\inf_{u\in U}F(t,x,y,z,u,v),\ (t,x,y,z)\in[0,T]\times\mathbb{X}\times\mathbb{R}\times\mathbb{R}^d.$$
Obviously,   $G$ is Lipschitz in $(y,z)$, uniformly
with respect to $(t,x)$ and there exists a constant $C>0$  such that
$
|G(t,x,y,z)|\leq C(1+\|x\|_t+|y|+|z|).
$ Since $F$ is continuous on the compact space $U\times V$, there exist two measurable mappings
$\bar{u}$ (resp. $\bar{v}):[0,T]\times\mathbb{X}\times\mathbb{R}\times
\mathbb{R}^d\rightarrow U$ (resp. $V$) such that
\begin{equation}\label{091201}
G(t,x,y,z)
=F(t,x,y,z,\bar{u}(t,x,y,z),\bar{v}(t,x,y,z)).
\end{equation}
Moreover, for all $(u,v)\in U\times V$, it holds
$$F(t,x,y,z,\bar{u}(t,x,y,z),v)\leq F(t,x,y,z,\bar{u}(t,x,y,z),\bar{v}(t,x,y,z))
\leq
F(t,x,y,z,u,\bar{v}(t,x,y,z)).$$
Then using  comparison theorem of RBSDEs and Theorem \ref{th6.1}, we get the results as follows.
\begin{theorem}\label{th6.3}
Suppose that the Isaacs condition holds. Then the value of the game problem \eqref{090801} exists, which can be characterized as follows
\begin{equation}\nonumber
\sup_{v\in\mathcal{V}_{\mathbb{F}}}
\inf_{u\in\mathcal{U}_{\mathbb{F}}}
J(u,v)=
\inf_{u\in\mathcal{U}_{\mathbb{F}}}
\sup_{v\in\mathcal{V}_{\mathbb{F}}}J(u,v)=
{P}_0,
\end{equation}
where $({P},{Q},A)$ is the unique solution of the following RBSDE
\begin{equation}\nonumber
\left\{
\begin{aligned}
&{P}_t= \Phi(X)+\int_t^TG(s,X,{P}_s,{Q}_s)ds +A_T-A_t-\int_t^T{Q}_sdW_s,\ t\in[0,T],\\
&{P}_t\geq h(t,X),\ \forall\ t\in[0,T],\ \text{a.s.},\ \int_0^T[{P}_{t}-h(t,X)]dA_t=0.
\end{aligned}
\right.
\end{equation}
Moreover, the admissible control pair $(u^*,v^*)\in\mathcal{U}_{\mathbb{F}}
\times\mathcal{V}_{\mathbb{F}}$  given as
$$u^*_t=\bar{u}(t,X,{P}_t,{Q}_t),\ v^*_t=\bar{v}(t,X,{P}_t,{Q}_t),
$$
forms a saddle point,
where the function $(\bar{u},\bar{v})$ is given in \eqref{091201}.
\end{theorem}

\bibliographystyle{SIAM}
\bibliography{mybib}

\end{document}